\documentclass[11 pt,reqno]{amsart}

\usepackage{graphics}
\usepackage{amsmath,amsthm,amssymb,amscd,graphicx}
\usepackage{color}
\usepackage{hyperref}
\usepackage{enumerate}
\usepackage{mathrsfs}
\usepackage{enumitem}
\usepackage{multirow}

\newtheorem{thm}{Theorem}[section]

\newtheorem{conjecture}[thm]{Conjecture}

 \theoremstyle{definition}
 
\newtheorem{definition}[thm]{Definition}
 \theoremstyle{remark}
 
 \newtheorem{remark}[thm]{Remark}
 
 \numberwithin{equation}{section}

\def\be#1 {\begin{equation} \label{#1}}
\newcommand{\ee}{\end{equation}}

\newcommand{\Real}{\mathbb R}
\newcommand{\Complex}{\mathbb C}

\newcommand{\eps}{\epsilon}
\newcommand{\Torus}{\mathbb T}

\newcommand{\brak}[1]{\langle #1 \rangle}
\newcommand{\grad}{\nabla}
\newcommand{\norm}[1]{\left\lVert #1 \right\rVert}
\newcommand{\abs}[1]{\left\vert#1\right\vert}
\newcommand{\set}[1]{\left\{#1\right\}}

\def\G{\mathcal G}
\def\cG{\mathcal G}

\begin{document}
\author{Jacob Bedrossian}
\address{4176 Campus Drive - William E. Kirwan Hall, College Park, MD 20742, USA}
\curraddr{}
\email{jacob@cscamm.umd.edu}

\author{Pierre Germain}
\address{∗Courant Institute of Mathematical Sciences, 251 Mercer Street, New York, NY 10012, USA}
\curraddr{}
\email{pgermain@cims.nyu.edu}

\author{Nader Masmoudi}
\address{‡Department of mathematics, New York University in Abu Dhabi, Saadyiat Island, Abu Dhabi, UAE}
\curraddr{}
\email{masmoudi@cims.nyu.edu}

\title{Stability of the Couette flow at high Reynolds numbers in 2D and 3D} 



\begin{abstract} 
We review works on the asymptotic stability of the Couette flow. 
The majority of the paper is aimed towards a wide range of applied mathematicians with an additional section aimed towards experts in the mathematical analysis of PDEs.   
\end{abstract}

\maketitle
\setcounter{tocdepth}{1}

\begin{quote}
\footnotesize\tableofcontents
\end{quote}

Here we provide a review on the stability of the 2D and 3D Couette flow in the Navier-Stokes equations at high Reynolds number, focusing mainly on the recent sequence of works \cite{BM13,BMV14,BGM15I,BGM15II,BGM15III,BVW16}. 
Our goal is to provide a general discussion of the wider physical context for these works and those which are related, and specifically, to discuss how they fit into the wider field of hydrodynamic stability. 
To this purpose,  we have attempted to make sections \S\ref{sec:Hydrostab} - \ref{sec:Statements} accessible to a wide range of applied mathematicians and mathematical analysts, whereas \S\ref{sec:Rigor} is targeted towards experts in the analysis of PDEs. 

\section{Hydrodynamic stability at high Reynolds number} \label{sec:Hydrostab}

\subsection{The Navier-Stokes and Euler equations}

We start by recalling the \textit{incompressible Navier-Stokes equations} in dimension $d=2$ or $3$, with inverse Reynolds number $\nu := \textbf{Re}^{-1} \geq 0$, 
\begin{equation}
\label{NS}
\tag{NS}
\left\{
\begin{array}{l}
\partial_t v + v \cdot \nabla v = - \nabla p + \nu \Delta v\\
\nabla \cdot v = 0 \\
v(t=0) = v_{in},
\end{array}
\right.
\end{equation}
where the velocity $v=\begin{pmatrix} v^1 \\ v^2 \end{pmatrix}$ or $v = \begin{pmatrix} v^1 \\ v^2 \\ v^3 \end{pmatrix}$takes values in $\mathbb{R}^d$, the pressure $p$ is a scalar, and $v\cdot \nabla$ stands for $\sum_{i = 1}^d v^i \partial_i$.
If the Reynolds number is infinity (hence $\nu = 0$), then the system is known as the \textit{incompressible Euler equations}, which reads:
\begin{equation}
\label{E}
\tag{E}
\left\{
\begin{array}{l}
\partial_t v + v \cdot \nabla v = - \nabla p \\
\nabla \cdot v = 0 \\
v(t=0) = v_{in},
\end{array}
\right.
\end{equation}
Both $v$ and $p$ are functions of the time and space variables $(t,x,y)$ or $(t,x,y,z)$, in space dimensions $2$ or $3$ respectively. 
Boundary conditions and external forces can also  be added to these equations, and are necessary for most non-trivial hydrodynamic stability questions. 

This article is concerned with hydrodynamic stability at high Reynolds number (i.e. in the singular limit $\nu \rightarrow 0$). 
See the texts \cite{DrazinReid81,SchmidHenningson2001,Yaglom12} for reviews on the wider theory of hydrodynamic stability; see especially \cite{Yaglom12} for an extensive review of the literature and a detailed account of the development of the field.
Here we will only present a brief introduction especially targeting mathematicians. 
The basic problem is to consider a given equilibrium\footnote{As remarked above, in order to get non-trivial equilibria, we usually need to impose boundary conditions or an external force field. However, let us ignore this for the moment, as it is not directly relevant to our discussions.}  $u_{E}$ for either \eqref{NS} or \eqref{E}, and to study the dynamics of solutions which are close to $u_E$ in a suitable sense (we will be more precise later). Especially, to answer the simple question of whether or not the flow remains close to $u_E$ in certain norms or not, and if so, which norms.
Hence, if we write the solution $v = u_{E} + u$, the exact evolution equations for the perturbation become 
\begin{align} \label{eq:perteqns} 
\left\{
\begin{array}{l}
\partial_t u + u_E \cdot \grad u + u \cdot \grad u_E + u \cdot \nabla u = -\nabla p + \nu \Delta u \\
\nabla \cdot u = 0 \\
u(t=0) = u_{in},
\end{array}
\right.
\end{align}
where $u_{in}$ is the initial perturbation. 
For hydrodynamic stability questions, naturally $u_{in}$ is assumed initially small in certain norms (in the limit $\nu \rightarrow 0$, the question can be sensitive to this choice, as we will see). 
The vast majority of work is focused on \emph{laminar} equilibria -- simple equilibria in which the fluid is moving in well-ordered layers (as opposed to equilibria with chaotic streamlines). 
Typical examples include shear flows in various geometries such as pipes and channels, vortices or vortex columns, and flows in concentric cylinders \cite{DrazinReid81,SchmidHenningson2001,Yaglom12}. 
However, even for these simple configurations, surprisingly little is understood about the near-equilibrium dynamics in the limit $\nu \rightarrow 0$, especially at the mathematically rigorous level.  
When studying this limit, it is natural to first consider $\nu = 0$, which, as a general rule, will be strictly easier than studying the singular limit. 
Indeed, significantly more is known about stability theory for the Euler equations ($\nu = 0$) than for Navier-Stokes $\nu > 0$ in the limit $\nu \rightarrow 0$. 

\subsection{Linear theory and its limitations} 
It is often implicitly presumed that linear stability is equivalent to nonlinear stability in practical settings for physically relevant systems.  
However, even the very early linear studies of Rayleigh \cite{Rayleigh80} and Kelvin \cite{Kelvin87} seemed in contradiction with experimental observations, particularly the famous experiments of Reynolds \cite{Reynolds83}. 
Indeed, in \cite{Reynolds83}, Reynolds pumped fluid through a pipe under various conditions and demonstrated that at sufficiently high Reynolds number, the laminar equilibrium becomes spontaneously unstable. 
Despite the observed nonlinear instability, to this day, there is no evidence of any linear instabilities in 3D pipe flow for any finite Reynolds number \cite{DrazinReid81,Yaglom12}.
Perhaps even more troubling is that even for laminar flows for which the linearization has unstable eigenvalues at high Reynolds numbers, experiments and computer simulations normally display instabilities which are \emph{different} than those predicted by the linear theory (and at lower Reynolds numbers) \cite{Carlson1982,Alavyoon1986,Yaglom12,SchmidHenningson2001}. 
These phenomena are known in fluid mechanics as \emph{subcritical transition} or \emph{by-pass transition}, and are completely ubiquitous in 3D hydrodynamics. 

%In this section we discuss the approach of studying the linearization of \eqref{eq:perteqns}. 
To understand the subtleties of subcritical transition, it is important to first be more precise about the meaning of `linear stability' and `nonlinear stability'. 
%The most classical technique for studying the stability of an equilibrium is to linearize the equations and study the associated linear problem, 
The linearization of \eqref{eq:perteqns} is obtained simply by dropping the quadratic nonlinearity in \eqref{eq:perteqns}:  
\begin{align}\label{eq:linpert} 
\left\{
\begin{array}{l}
\partial_t u + u_E \cdot \grad u + u \cdot \grad u_E = -\nabla p + \nu \Delta u \\
\nabla \cdot u = 0 \\
u(t=0) = u_{in}.
\end{array}
\right.
\end{align}
Write \eqref{eq:linpert} abstractly as 
\begin{align}
\left\{
\begin{array}{l}
\partial_tu = \mathcal{L}_E u \\ 
u(t=0) = u_{in}.
\end{array}
\right.
\end{align}
The simplest notion of linear stability is \emph{spectral stability:}

\medskip
\noindent
\fbox{ \parbox{.98\linewidth}{
\begin{definition}[Spectral stability] \label{def:spec}
Given a Hilbert space $H$, let $\sigma(\mathcal{L}_E)$ be the spectrum of $\mathcal{L}_E$. 
The equilibrium $u_E$ is called \emph{spectrally stable} if $\sigma(\mathcal{L}_E) \cap \set{c \in \Complex: \textup{Re}\, c > 0} = \emptyset$. 
The equilibrium is called \emph{spectrally unstable} if $\sigma(\mathcal{L}_E) \cap \set{c \in \Complex: \textup{Re}\, c > 0} \neq \emptyset$ and is called \emph{neutrally stable} if it is spectrally stable but $\sigma(\mathcal{L}_E) \cap \set{c \in \Complex: \textup{Re}\, c = 0} \neq \emptyset$. 
\end{definition}}}
\medskip

Many classical theories of hydrodynamic stability are focused on studies of spectral stability (see \cite{Yaglom12} and the references therein). 
In the context of shear flows at high Reynolds number, the most famous results are Rayleigh's inflection point theorem \cite{Rayleigh80,DrazinReid81} and its refinement due to Fj{\o}rtoft \cite{Fjortoft1950}, regarding the stability of shear flows in the 2D Euler equations ($\nu = 0$). 
These results (unfortunately) extend to planar 3D shear flows via Squire's theorem \cite{Squire1933,DrazinHoward1966,DrazinReid81}. 
Squire's theorem for inviscid flows states: If the 3D Euler equations linearized around the shear flow $u_E = (f(y),0,0)$ have unstable eigenvalues, then so do the 2D Euler equations linearized around $u_E = (f(y),0)$ (the converse implication being obvious).  
Therefore, for any shear flow of this type, spectral stability in 2D implies spectral stability in 3D. 
This is `unfortunate' because it gives the false impression that most of the interesting aspects of hydrodynamic stability can be found in the 2D equations. 
As we will see, hydrodynamic stability in 2D and 3D, at high Reynolds numbers, are quite distinct. 

Spectral stability is not the only relevant definition of linear stability. 
In particular, we have the following distinct definition, which was suggested as more natural in fluid mechanics by Kelvin \cite{Kelvin87}, Orr \cite{Orr07}, and later Case \cite{Case1960}and Dikii \cite{Dikii1960}.  

\medskip
\noindent
\fbox{ \parbox{.98\linewidth}{
\begin{definition}[Lyapunov linear stability] \label{def:linstab}
Given two norms $X$ and $Y$ (often assumed the same), the equilibrium $u_E$ is called \emph{linearly stable} (from $X$ to $Y$) if
\begin{align*}
\norm{u}_{Y} \lesssim \norm{u_{in}}_X,
\end{align*}
where $u$ solves \eqref{eq:linpert}. 
\end{definition} }}
\medskip

Even for finite-dimensional systems of ODEs, Definitions \ref{def:spec} and \ref{def:linstab} are not equivalent for neutrally stable, non-diagonalizable linear operators (the infinite dimensional analogue being non-normal operators).
The canonical example of this phenomenon, which is directly relevant for hydrodynamic stability as we will see, is given by the following linear system:
\begin{align*}
\partial_t \begin{pmatrix} x_1 \\ x_2 \end{pmatrix} = \begin{pmatrix} 0  & 1 \\ 0 & 0 \end{pmatrix} \begin{pmatrix} x_1 \\ x_2 \end{pmatrix}. 
\end{align*}
The matrix $\begin{pmatrix} 0  & 1 \\ 0 & 0 \end{pmatrix}$ is not diagonalizable. It has spectrum reduced to $\{ 0 \}$, and is therefore neutrally stable; but it is not Lyapunov stable: solutions are given by 
\begin{align*}
\left\{ \begin{array}{l} x_1(t) = x_1(0) + t x_2(0) \\ x_2(t) = x_2(0), \end{array} \right.
\end{align*}
and therefore grow.
The above example can be modified to make it spectrally stable and Lyapunov linearly stable: 
\begin{align}
\partial_t \begin{pmatrix} x_1 \\ x_2 \end{pmatrix} = \begin{pmatrix} -\nu  & 1 \\ 0 & -\nu \end{pmatrix} \begin{pmatrix} x_1 \\ x_2 \end{pmatrix}. \label{def:transgrwth}
\end{align}
As pointed by Orr in 1907 \cite{Orr07}, non-normality implies that a large transient growth is possible. As $\nu \to 0$, the system degenerates, leading to the optimal estimate
\begin{align*}
\sup_t \norm{(x_1(t),x_2(t))} \lesssim \nu^{-1}\norm{(x_1(0),x_2(0))}. 
\end{align*}
Hence, while the problem is linearly stable for all $\nu > 0$, the solutions are (in general) growing very large as $\nu \rightarrow 0$. 
As Orr pointed out in \cite{Orr07}, we could expect that any prediction of nonlinear stability from a linear system of this type should be problematic as $\nu \rightarrow 0$, a fact which has been greatly expanded on in modern times (see below). 
Many problems of hydrodynamic stability display transient behavior reminiscent of \eqref{def:transgrwth} in the limit $\nu \rightarrow 0$. 
See below and \cite{TTRD93,SchmidHenningson2001,Trefethen2005,PradeepHussain06} and the references therein for more details.   

 \subsection{Quantitative nonlinear stability and `transition thresholds'} \label{sqs}
The limitations of linear theory being increasingly clear, we must turn to more nonlinear stability theories in order to connect with physical observations. 
The simplest notion of nonlinear stability is that of Lyapunov.

\medskip
\noindent
\fbox{ \parbox{.98\linewidth}{
\begin{definition}[Lyapunov (nonlinear) stability] \label{def:nonsta0}
Given two norms $X$ and $Y$ (often assumed the same), the equilibrium $u_E$ is called \emph{stable} (from $X$ to $Y$) if for all $\epsilon > 0$, there exists a $\delta > 0$ such that if 
\begin{align*}
\norm{u_{in}}_X < \delta \quad \implies \quad \norm{u(t)}_Y < \epsilon \quad \mbox{for all}\; t>0.
\end{align*} 
We say $u_E$ is \emph{unstable} if it is not stable. 
\end{definition}  }}
\medskip

There are many  works in fluid mechanics and related fields dedicated to proving spectral instability and to proving that spectral instability implies nonlinear instability, see e.g. \cite{FriedlanderVishik1991,FriedlanderStraussVishik1997,ShvydkoyFriedlander2005,FriedlanderEtAl2006,Grenier2016,Grenier00} to name a few in the mathematics community. This theory has led to the clearer understanding of many physical phenomena in hydrodynamics \cite{Rayleigh1916,Heisenberg1984,lin1944,DrazinReid81} but it does not help us to understand the phenomenon of subcritical transition as it does not help us explain instabilities which are not directly associated with a spectral instability. 

Probably the most flexible and powerful general theory for nonlinear stability results in hydrodynamics is that of Arnold for the 2D Euler equations, introduced in \cite{Arnold1965}, now usually called the Energy-Casimir method.
This theory provides nonlinear stability in the sense of Definition \ref{def:nonsta} for $X = Y = H^1$ for certain equilibria, such as shear flows satisfying certain convexity hypotheses and vortices \cite{WanPulvirenti85,Dritschel88,MarchioroPulvirenti12}.
Variations have been extended to a great many other settings, for example magneto-hydrodynamics and kinetic theory; see the review article \cite{HolmEtAl85} and the more recent works in stellar mechanics \cite{Guo1999,Guo2001} for more information. 
The fundamental idea behind this technique is to use the large number of extra conservation laws available to 2D Euler, the so-called \emph{Casimirs}, to find a conserved quantity which is locally convex, coercive in $H^1$, and for which the equilibrium is the local minimizer. 
There are several limitations of the theory: it is mostly restricted to two dimensional settings and infinite Reynolds number (recall that, in the presence of boundaries, one can find examples of 2D shear flows in a channel which are stable at infinite Reynolds numbers by Arnold's theory, but are spectrally and nonlinearly unstable in $L^2$ at arbitrarily high (but finite) Reynolds numbers \cite{Grenier2016,FriedlanderEtAl2006}). 

The above discussion suggests that we are still lacking a good mathematical understanding of subcritical transition (recall that, for equilibria which are linearly stable, subcritical transition refers to the spontaneous transition to a turbulent state often observed at high Reynolds numbers in experiments and computer simulations). 
It has also been observed, even by Reynolds \cite{Reynolds83}, that at which Reynolds number and precisely how it occurred, seems unusually sensitive to the details of the experimental set-up. In \cite{Kelvin87}, Kelvin suggested the idea that maybe such equilibria are nonlinearly stable for all finite Reynolds number in the sense of Definition \ref{def:nonsta0}, but that the system becomes increasingly sensitive to perturbations as Reynolds number increased.
This suggests that what we should study is a quantification of $\delta$ in Definition \ref{def:nonsta0}. 

\medskip
\noindent
\fbox{ \parbox{.98\linewidth}{
\begin{definition}[Quantitative asymptotic stability] \label{def:nonsta}
Given two norms $X$ and $Y$, the equilibrium $u_E$ is called \emph{asymptotically stable} (from $X$ to $Y$) with exponent $\gamma$ if
\begin{align*}
\norm{u_{in}}_X \ll \nu^\gamma \quad \implies \quad 
\left\{ \begin{array}{ll}
\| u(t) \|_Y \ll 1 & \mbox{for all $t>0$} \\
\| u(t) \|_Y \to 0  & \mbox{as $t \to \infty$}.
\end{array} \right.
\end{align*} 
We say $u_E$ is \emph{unstable} if it is not stable. 
\end{definition}  }}
\medskip

We remark that the polynomial dependence $\nu^\gamma$, as well as the necessity of two norms $\| \cdot \|_X$ and $\| \cdot \|_Y$ is dictated in part by empirical observation.
As we will see, the choice of `final norms' $\|\cdot\|_Y$ is severely constrained by the linear dynamics, however, the initial norms $\|\cdot\|_X$ are not. 
In fact, there is no unique, natural physical choice of $X$.
Moreover, numerical experiments \cite{ReddySchmidEtAl98} and weakly nonlinear heuristics suggest that $\gamma$ should depend non-trivially on the choice of norms. 
While this has not yet been proved for the Navier-Stokes equations, a similar sensitivity to initial regularity was recently established for the Landau damping of Vlasov-Poisson in a periodic box \cite{Bedrossian16}; see \S\ref{sec:LandauMixing} below for more discussion on how Landau damping is related to hydrodynamic stability questions.

In case an equilibrium is asymptotically stable, the above definition raises the following question (given $X$ and $Y$):
\begin{equation}
\label{question}
\tag{Q1}
\mbox{\textbf{What is the smallest exponent $\gamma>0$ for which $u_E$ is quantitative asymptotic stable?}}
\end{equation}
As suggested above, the motivation for Definition \eqref{def:nonsta} and Question \eqref{question} go all the way back to Kelvin; however, modern authors have expanded on these concepts further: see e.g. \cite{TTRD93,ReddySchmidEtAl98,SchmidHenningson2001} and the references therein. 
Another very important question, which is a little more vague, is to understand precisely how the instability occurs near the stability threshold. 
\begin{align}
\label{question2}
\tag{Q2}
& \mbox{\textbf{If $\gamma$ is the optimal exponent in~\eqref{question},}} \\
& \mbox{\textbf{\quad for a $\gamma' < \gamma$, classify all solutions satisfying}} \nonumber\\
& \qquad   \qquad \qquad \nu^\gamma \lesssim \| u_{in} \|_X \lesssim \nu^{\gamma'} \quad \mbox{and} \quad \sup_{t \geq 0}\| u(t) \|_Y  \gtrsim 1. 
\end{align}

When considering only the question of spectral stability or instability of planar shear flows, Squire's theorem suggests that one can focus exclusively on the 2D case (see e.g. \cite{DrazinReid81}). 
However, one gets very different answers in 2D and 3D when considering the nonlinear questions \eqref{question} and \eqref{question2}. 
Indeed, it is often observed that subcritical instabilities appear at lower Reynolds numbers than eigenvalue instabilities, and that even when 2D eigenvalue instabilities are present, fully 3D, non-modal instabilities are often still observed to be dominant \cite{SchmidHenningson2001,Yaglom12}. 
Similarly, we see major differences when comparing our 2D and 3D works on planar Couette flow for sufficiently smooth data (see below and \cite{BM13,BMV14} for the 2D works and \cite{BGM15I,BGM15II} for 3D). 
One cannot escape the reality that \emph{any study of 2D hydrodynamic stability at high Reynolds number is of very limited physical relevance to 3D flows}.

\subsection{Inviscid damping and asymptotic stability at infinite Reynolds number}

While question~\eqref{question} explicitly concerns only $\nu>0$, it is natural to ask, if $\gamma = 0$, whether some kind of stability holds for the inviscid problem $\nu = 0$. 
For non-trivial shear flows, one can verify that this is not possible unless $d = 2$ (see \S\ref{PED3} for more discussion); indeed, it is reasonable to imagine that every non-trivial, laminar 3D equilibrium is nonlinearly unstable in the 3D Euler equations (regardless of how strong one is willing to choose the initial norm $X$). 
However, for $d=2$, it turns out in some cases to be possible not only to have Lyapunov stability, but also to identify the long-term dynamics and have a kind of \emph{asymptotic} stability even if $\nu = 0$.
In \cite{BM13}, it was demonstrated that near the shear flow $u_E = (y,0)$ (with $(x,y) \in \mathbb T \times \Real$), for $u_{in}$  sufficiently small and sufficiently smooth, there exists some $u_\infty$ such that
\begin{align*}
u(t,x,y) \rightarrow \begin{pmatrix} u_{\infty}(y)\\ 0 \end{pmatrix} \quad \mbox{as $t \to \infty$}
\end{align*}
\emph{strongly} in $L^2$ and weakly in $H^1$ (and analogously also $\lim_{t \rightarrow -\infty} u(t) = (u_{-\infty}(y),0)$ for some $u_{-\infty}$).
Notice that this is \emph{much} stronger than the kinds of stability results one derives via the energy-Casimir method. 
This convergence back to equilibrium, despite time-reversibility and the lack of dissipative mechanisms, is known as \emph{inviscid damping} and is a close relative of Landau damping in plasma physics. It was proved that Landau damping provides a similar stability for Vlasov-Poisson in $\Torus^d$ in Mouhot and Villani's breakthrough work \cite{MouhotVillani11}; (see \S\ref{sec:LandauMixing} for more discussion on Landau damping).
Inviscid damping was observed in the linearized Euler equations first by Orr in 1907 \cite{Orr07}; see \S\ref{PED2} for more discussion on this effect. 
Landau damping is considered fundamental to understanding the dynamics of collisionless and weakly collisional plasmas by the physics community \cite{Ryutov99,BoydSanderson} and it has been speculated that inviscid damping should play a related role in understanding the dynamics of the 2D Euler equations \cite{GilbertBassom98,Gilbert88,MelanderEtAl1987,YaoZabusky96,SchecterEtAl00}, with applications to cyclotron dynamics \cite{CerfonEtAl13} and atmospheric sciences \cite{McCalpin87,SmithMontgomery95,Wang2008}. 
See \cite{LinZeng11,WeiZhangZhao15,zillinger2016,Zillinger2016circ,WeiZhangZhao2017,LinXu2017,YangLin16,WZZK2017,BCZV17} for mathematical works on inviscid damping in the linearized 2D Euler equations. 

%Lastly, we note that it should not be surprising that the $\nu > 0$ case studied in \cite{BMV14} is actually harder than the $\nu = 0$ case studied previously in \cite{BM13} -- the former is a `long-time inviscid limit' of the latter. 

\subsection{The Couette flow}
It turns out that the simplest equilibrium to study is the plane, periodic Couette flow: 
\begin{equation*}
u_E(x,y) = 
\begin{pmatrix} y \\ 0 
\end{pmatrix}\;
\;\;(\mbox{if}\; d=2)
\qquad \quad
u_E(x,y,z) = 
\begin{pmatrix} y \\ 0 \\ 0
\end{pmatrix}\;\;(\mbox{if}\;d=3).
\end{equation*}
Here `plane' differentiates this flow from the Taylor-Couette flow, the analogous equilibrium between concentric, rotating cylinders \cite{DrazinReid81}, while `periodic' refers to the choice $x \in \mathbb{T}$ (the circle, which we identify with $[0,2\pi]$ with periodic boundary conditions). 
Specifically, we consider the following geometry for the domain:
\begin{itemize}
\item In dimension $d=2$, $(x,y) \in \mathbb{T} \times \mathbb{R}$.
\item In dimension $d=3$, $(x,y,z) \in \mathbb{T} \times \mathbb{R} \times \mathbb{T}$.
\end{itemize}
>From now on we will simply refer to this flow as the Couette flow. 

Why study the Couette flow in this domain specifically? There are several good reasons: 
\begin{itemize}
\item The Couette flow is a canonical problem and has played the role of a benchmark for different approaches to understanding hydrodynamic stability for well over a century (see \cite{DrazinReid81,Yaglom12,SchmidHenningson2001} and the references therein). This is due in part to its simplicity, and to the fact that it is a solution of the Navier-Stokes equations for all $\nu \geq 0$. 
\item The domain of infinite extent in the $y$ direction avoids the presence of walls, which can add very subtle complications to hydrodynamic stability problems (see also \S\ref{sec:openboundariesx} for more discussion). 
\item The shear gives rise to a mixing phenomenon (sometimes referred to as `filamentation') which has profound implications for the dynamics and is mostly a stabilizing effect; however, this effect is very weak at low wavenumbers in $x$. The periodic domain in the $x$ direction removes such problematic wave numbers. 
The periodicity is most similar to the mixing which occurs around a radially symmetric vortex, however, studying the stability of vortices at high Reynolds number is significantly more difficult (see \cite{BCZV17,LiWeiZhang2017,Gallay2017}). 
\item The reduction to periodic boundary conditions in $x$, and infinite extent in $y$ is common in the physics literature both for formal asymptotic analysis~\cite{TTRD93,Chapman02,SchmidHenningson2001} and (at least the periodicity in $x$) for numerical studies. 
\end{itemize}

\subsection{Summary of the known stability results on the plane, periodic Couette flow}

Many works have studied the stability of the plane, periodic Couette flow and its variants in attempt to get answers to Questions \eqref{question} and \eqref{question2}. 
For experiments on Couette flow, see \cite{Tillmark92,Daviaud92,bottin98}, and for computational studies, see e.g. \cite{Orszag80,ReddySchmidEtAl98,DuguetEtAl2010}. 
In \cite{TTRD93}, the authors suggested a weakly nonlinear toy model which couples transient growth in the linear problem to a caricature of the nonlinearity in attempt to capture some of the aspects of subcritical transition. For more work along this general line, see e.g.  \cite{Gebhardt1994,BDT95,Waleffe95,BT97,LHRS94,Chapman02}. 
As far as rigorous analysis goes, there is the work of Romanov \cite{Romanov73} who first proved nonlinear stability at all Reynolds number for the Couette flow in an infinite pipe geometry. Later, the works of \cite{KreissEtAl94,Liefvendahl2002} estimate $\gamma$ rigorously assuming a numerically obtained resolvent estimates.  
See \cite{Yaglom12} for further references on Couette flow as well as references focusing on related flows, such as the plane Poiseuille flow and the cylindrical pipe flow.

The table \ref{fig1} recapitulates known answers to Question \eqref{question} as deduced in the sequence of works \cite{BM13,BMV14,BGM15I,BGM15II,BGM15III,BVW16}. 
Relevant norms we will use to measure the initial data are given by Sobolev norms and ``Gevrey norms'', denoted $H^s$ and $\mathcal{G}^{\lambda,s}$ respectively, defined via the Fourier transform (here $\mathbf \xi$ denotes the wavenumber; see \S\ref{sec:FourierConv})
\begin{align*}
\norm{f}_{H^s} & = \left(\sum \int \abs{\hat{f}(\mathbf \xi)}^2 (1 + \abs{\mathbf \xi})^{2s} d\mathbf \xi\right)^{1/2} \\ 
 \norm{f}_{\mathcal{G}^{\lambda;s}} & = \left(\sum \int \abs{\hat{f}(\mathbf \xi)}^2 e^{2\lambda \abs{\mathbf{\xi} }^s} d\mathbf \xi\right)^{1/2}. 
\end{align*}
Roughly speaking, a finite Gevrey $\mathcal{G}^{\lambda;s}$ norm requires an average decay slightly stronger than $ e^{-\lambda |\xi|^s}$ in frequency, while a finite Sobolev norm requires an average decay slightly stronger than $|\xi|^{-s}$.
Accordingly, perturbations which are small in a Gevrey norm are much smoother than those perturbations which are only small in a Sobolev norm. 
Below, we will broadly divide  results into those which measure the initial data with a Sobolev norm versus those results which measure the initial data with the much stronger Gevrey norms. More specific requirements on the exact norms are discussed further below in \S\ref{sec:Statements}. 

\medskip
\begin{center}
\begin{figure} 
\begin{tabular}{|c|c|c|c|} 
\hline
\multirow{2}{*}{dimension $d=2$} & Gevrey & $\gamma = 0$ & \cite{BM13,BMV14} \\
\cline{2-4}
& Sobolev & $\gamma \leq \frac{1}{2}$ & \cite{BVW16} \\
\hline
\multirow{2}{*}{dimension $d=3$} & Gevrey & $\gamma \leq 1$ & \cite{BGM15I,BGM15II} \\
\cline{2-4}
& Sobolev & $\gamma \leq \frac{3}{2}$ & \cite{BGM15III}\\
\hline
\end{tabular}
\caption{Progress on question \eqref{question} for the plane, periodic Couette flow with $y \in \Real$.} \label{fig1}
\end{figure} 
\end{center}

Physicists are in agreement that (for 3D) $\gamma \leq 1$ and the only question they considered is whether $\gamma =1 $ or $\gamma > 1$, as conjectured by many (see \S\ref{sec:streaks} and \cite{TTRD93,ReddySchmidEtAl98,Chapman02,SchmidHenningson2001}). However, strictly speaking a full mathematical proof that $\gamma \geq 1$ is not currently carried out to our knowledge.  
Computer simulations of `smooth' initial perturbations have often returned values slightly over-estimating the (suspected) answer of $\gamma = 1$ for Gevrey spaces (see \cite{ReddySchmidEtAl98,DuguetEtAl2010} and the references therein). 
Indeed, the general trend has been a steady decrease down to $1$, which was only reported relatively recently in \cite{DuguetEtAl10}. 
However, $\gamma = 31/20$ was reported for `rough' initial perturbations in \cite{ReddySchmidEtAl98}; this exponent seems more robust to computation - and matches closely the value $\frac{3}{2}$ appearing in the above table. 
We are not currently aware of any 2D simulations of this type.  
In addition to the estimates on $\gamma$, these works \cite{BM13,BMV14,BVW16,BGM15I,BGM15II,BGM15III} also provide detailed information regarding the long-time dynamics of the solutions. 
In fact, understanding these dynamics is crucial to getting accurate estimates of $\gamma$. 

>From a mathematically rigorous perspective, we have significantly less information regarding the classification of instabilities~\eqref{question2}. 
In general, this problem appears to be significantly harder than \eqref{question}.
 The only (rigorous) results in this direction are that of \cite{BGM15II}, which consider
the Gevrey norm in 3D with initial data $\nu \lesssim \norm{u_{in}}_{\mathcal{G}^{s,\lambda}} \lesssim \nu^{2/3+\delta}$ for $\delta > 0$ small. 
In this work, we prove that only one instability is possible. The instability is driven by a distinctly 3D, linear, non-normal transient growth known as the lift-up effect \cite{EllingsenPalm75,landahl80}; see \S\ref{PED3} below for more discussion. 
A related work is that of \cite{LinZeng11}, in which the authors prove that inviscid damping matching the linearized problem is false in general for the nonlinear 2D Euler equations with perturbations of Couette flow arbitrarily small in $H^s$ for $s < 3/2$.

\begin{remark}
It is inaccurate to suggest that all physical settings should be compared to Gevrey or Sobolev (or $L^2$) norms.
One should instead measure the size of the perturbations in the various norms after restricting to scales larger than a suitable dissipative length-scale. 
%the Fourier transform of the perturbations down to a suitable dissipative length-scale and compare the size of the data restricted to scales larger than the dissipative scales in the differents topologies.
Essentially, below this scale, viscosity dominates and the difference between the norms is no longer relevant. 
%Unfortunately, this kind of analysis seems rather impractical to carry out for a real system. 
%However, 
\end{remark}

\subsection{Notations and Fourier analysis conventions} \label{sec:FourierConv}
We will use the Fourier transform frequently in 2D and 3D.  
In 2D we use 
\begin{align*}
\widehat{f}(k,\eta) = \frac{1}{2\pi}\int_{\mathbb{T} \times \mathbb{R}} e^{i(kX + \eta Y )} f(X,Y,Z) \,dX\,dY.
\end{align*}
and in 3D we use 
\begin{align*}
\widehat{f}(k,\eta,\ell) = \frac{1}{(2\pi)^{3/2}}\int_{\mathbb{T} \times \mathbb{R} \times \mathbb{T}} e^{i(kX + \eta Y + \ell Z )} f(X,Y,Z) \,dX\,dY \, dZ.
\end{align*}
Note that $k,\ell \in \mathbb{Z}$ whereas $\eta \in \Real$. 
It will be very important in the following to distinguish between zero and non-zero frequencies in $x$ and $X$ ($k=0$ or $k \neq 0$). We will denote the orthogonal projection onto the modes $f_0$ and $f_{\neq}$ respectively:
\begin{align}
f_0 = \frac{1}{2\pi}\int_{\mathbb{T}} f\,dX \quad \mbox{and} \quad f_{\neq} = f - f_0.\label{def:prj}
\end{align}

\section{Linear dynamics in dimension $d=2$}

\label{PED2}
In two dimensions, it is most natural to consider the problem in the vorticity formulation. Recall that the vorticity $\omega$ of the perturbation is given by
$$
\omega = \partial_x u^2 - \partial_y u^1.
$$
We will also use the stream function $\phi$, which satisfies
\begin{equation}
\label{defphi}
\Delta \phi = \omega, \qquad u = \nabla^\perp \phi.
\end{equation}
The equation in vorticity form reads
\begin{align}
\label{def:vortform}
\left\{ \begin{array}{l}
\partial_t \omega + y \partial_x \omega + \nabla^\perp \phi \cdot \nabla \omega = \nu \Delta \omega \\
\Delta \phi = \omega.
\end{array} \right.
\end{align}

\subsection{Coordinate transform} 
\label{NC2D}
For small data and small $\nu$, we expect the transport by the Couette flow in \eqref{def:vortform} to be dominant. 
Hence, it is natural to switch to coordinates in which this transport is canceled or `modded out'.
To that effect, set 
\begin{equation} \label{2dcoords}
\left\{
\begin{array}{l}
X = x-ty \\
Y = y.
\end{array}
\right.
\end{equation}
We will follow the convention that dependent variables which are capitalized are taken in the new coordinates $(X,Y)$, whereas small font is reserved for the original coordinates $(x,y)$. Hence, for instance
\begin{subequations}
\begin{align}
U(t,X,Y) & = u(t,x,y) \\
\Omega(t,X,Y) & = \omega(t,x,y).
\end{align}
\end{subequations}
Under this change of coordinates, differential operators transform as follows
\begin{align*} 
\left\{
\begin{array}{l}
\nabla u = \nabla_L U \\
\Delta u = \Delta_L U
\end{array}
\right.
\quad \mbox{where} \quad
\nabla_L = \begin{pmatrix} \partial_X \\ \partial_Y - t \partial_X \\ \end{pmatrix}
\quad
\mbox{and}
\quad
\Delta_L = \nabla_L \cdot \grad_L = \partial_X^2 + (\partial_Y - t \partial_X)^2.
\end{align*}
In the new coordinates \eqref{2dcoords}, the equation satisfied by $\Omega$ reads
\begin{align}
\label{pinguin}
\left\{ \begin{array}{l}
\partial_t \Omega + \nabla^\perp \Phi \cdot \nabla \Omega  = \nu \Delta_L \Omega \\
\Delta_L \Phi = \Omega.
\end{array} \right.
\end{align}
Note that there was a crucial cancellation: the terms involving $\pm t\partial_X \Phi \partial_X \Omega$ cancel. Such a fortuitous cancellation will not hold in 3D. 
Linearizing this for small perturbations gives the very simple 
\begin{align} \label{linOmega2D} 
\partial_t \Omega & = \nu \Delta_L \Omega.
\end{align}

\subsection{Enhanced dissipation} \label{sec:ED2D}
The fact that the dissipation of a passive scalar can be enhanced via advection by an incompressible velocity field is a well-known effect, going sometimes by the name of `shear-diffuse mechanism', `relaxation enhancement', and `enhanced dissipation'. It has been studied previously in linear problems in both the mathematics literature \cite{CKRZ08,Zlatos2010,BeckWayne11,VukadinovicEtAl2015,BCZGH15,BCZ15,BZ09,IYM17,LiWeiZhang2017} (see also \cite{GallagherGallayNier2009}) and in the physics literature \cite{Lundgren82,RhinesYoung83,DubrulleNazarenko94,LatiniBernoff01,BernoffLingevitch94}.
We remark that Lundgren \cite{Lundgren82} considered possible implications for the behavior of vortex filaments in turbulent flow whereas Dubrulle and Nazarenko \cite{DubrulleNazarenko94}
suggested the possibility that the enhanced dissipation effect could have an important impact on Questions \ref{question} and \ref{question2}. 
This effect can be explained as follows: the advection transfers enstrophy to very high frequencies where it is then more rapidly dissipated.
 It was Kelvin \cite{Kelvin87} who first made the observation, when he solved \eqref{linOmega2D} via Fourier analysis in 1887.

First, the zero frequencies in $X$ are governed by the simple heat equation (recall that $\Omega_0 = \int \Omega\,dX$): 
$$
\partial_t \Omega_0 = \nu \partial_{YY} \Omega_0. 
$$
Note that $\Omega_0(t,y)$ defines the vorticity of a shear flow by the Biot-Savart law in \eqref{linOmega2D}.  The dissipation time-scale is $\sim \nu^{-1}$ and there is obviously no enhancement for these modes.

Next, we consider the non-zero frequencies in $X$. 
Applying the Fourier transform to \eqref{linOmega2D} gives: 

$$
\partial_t \widehat{\Omega}(k,\eta) = - \nu \left[ k^2 + (\eta - kt)^2  \right] \widehat{\Omega}(k,\eta). 
$$

If $k\neq 0$, then integrating the above ODE gives: 
\begin{align*}
\widehat{\Omega}(k,&\eta,\ell)  = \exp\left(-\nu \int_0^t (k^2 + (\eta-ks)^2) \,ds \right)\widehat{\Omega_{in}}(k,\eta) \\
& \implies |\widehat{\Omega}(k,\eta,\ell)| \lesssim \exp(-c \nu k^2 t^3)\abs{\widehat{\Omega_{in}}(k,\eta)},
\end{align*}
for a fixed, small constant $c>0$ ($1/12$ suffices). 
Using that $x \in \Torus$, this implies that, in any reasonable norm (by the Biot-Savart law $U = \grad^\perp_L \Delta_L^{-1}\Omega$), 
$$
\norm{\Omega_{\neq}(t)} + \| U_{\neq}(t) \| \lesssim e^{-c \nu t^3}\norm{\Omega_{in}}.
$$
The exponent $\nu t^3$ gives a dissipation time scale $\sim \nu^{-1/3}$. This is very short compared to the dissipation time scale $\nu^{-1}$ which is observed for zero frequencies in $X$. 
In more general linear advection-diffusion problems or linearized Navier-Stokes equations, precise, quantitative estimates on this `enhanced dissipation' effect are significantly harder to obtain see e.g. \cite{CKRZ08,BeckWayne11,BCZ15,LiWeiZhang2017,IYM17}. 
%The effect can  also give rise to hypoellipticity in Gevrey classes \cite{BCZ15}. 
Notice  that the effect becomes weak at low wavenumbers; see \S\ref{sec:openboundariesx} for more discussion.  

The mixing effect which gives rise to enhanced dissipation is a very common phenomenon in high Reynolds number fluid mechanics. 
One can see it when stirring milk into coffee. 

\begin{figure}[ht]
        \begin{center}$ 
\begin{array}{ll}
        \includegraphics[width= 0.20 \textwidth]{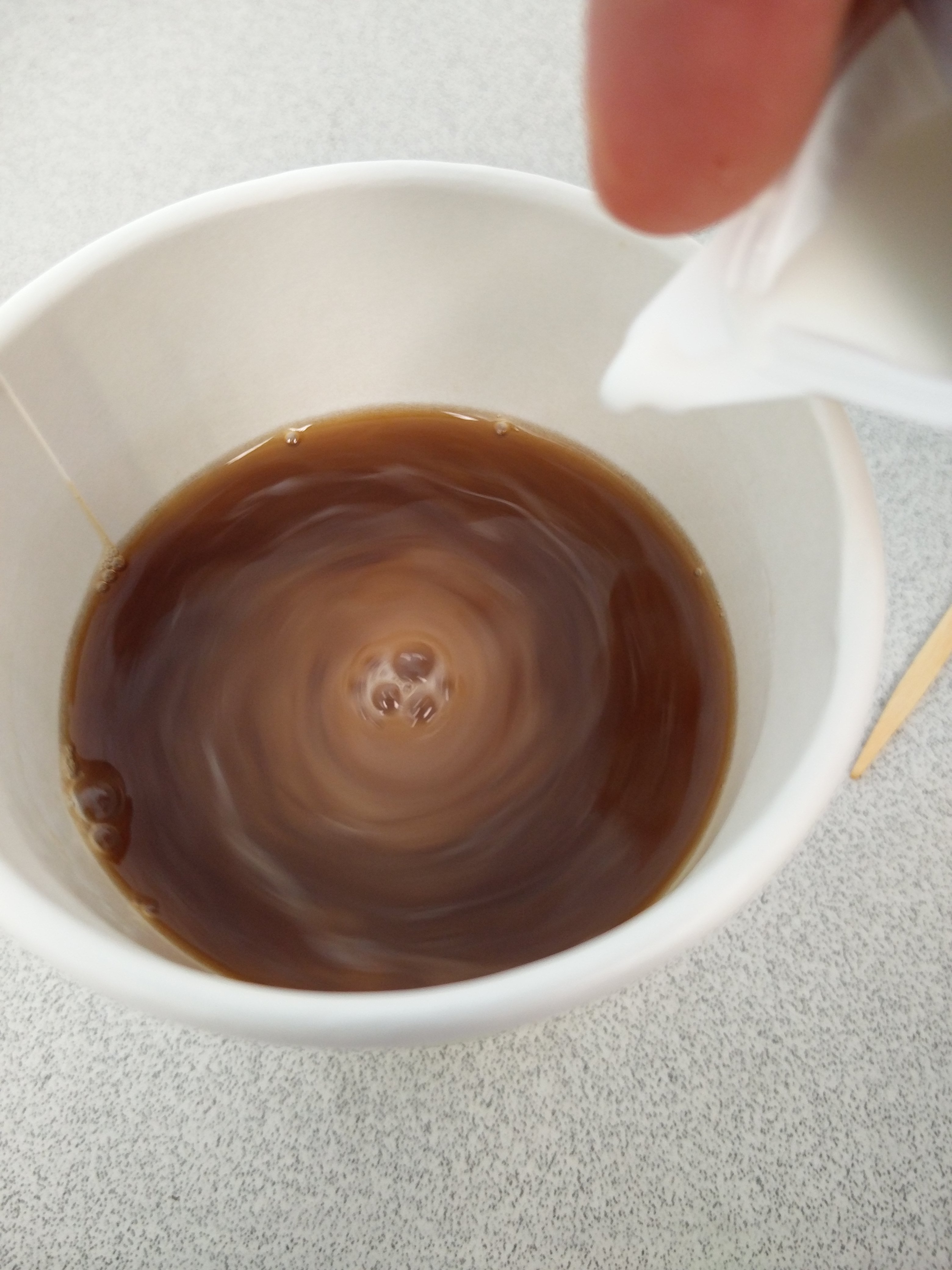} & %\hspace{2cm}
        \includegraphics[width= 0.20 \textwidth]{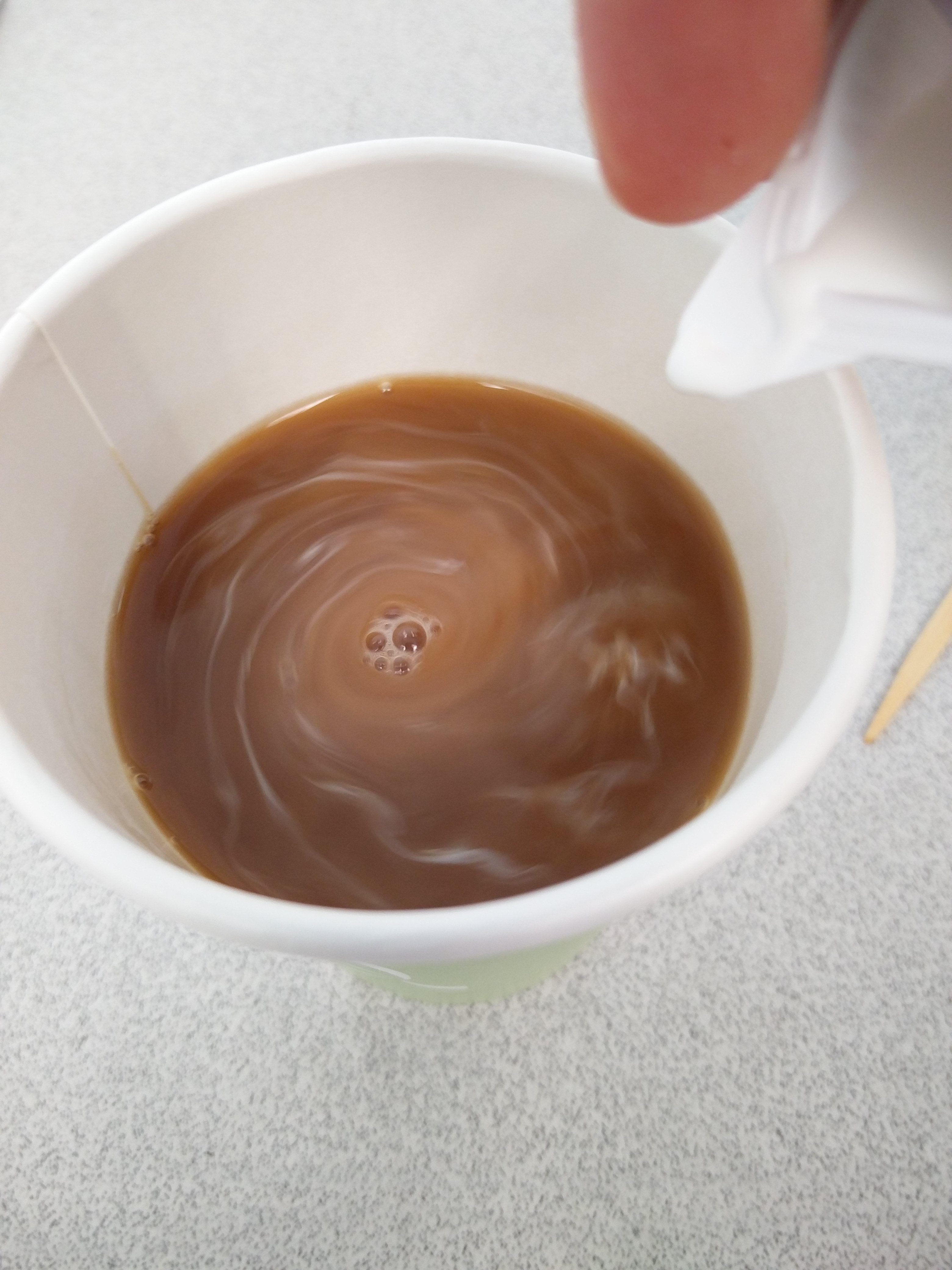}
        %\caption{Mmm coffee.} \label{fig:Coffee}
\end{array}
$\end{center}
\caption{Milk mixing into an approximately axisymmetric vortex.} \label{fig:Coffee}  
\end{figure}

\subsection{Inviscid damping and the Orr mechanism} \label{IDOM2D}
The enhanced dissipation already implies that the Couette flow should have a stabilizing effect. 
However, there is an additional, more subtle, inviscid effect. 
This was discovered  by Orr in 1907 \cite{Orr07} and is often known now as inviscid damping, in part due to its connection with Landau damping in the kinetic theory of plasmas (see \S\ref{sec:LandauMixing}).
>From the Biot-Savart law, we have the following: 
\begin{align*} 
\widehat{U}(t,k,\eta) = -i \begin{pmatrix}kt-\eta \\ k \end{pmatrix} \frac{1}{k^2 + (\eta-kt)^2} \widehat{\Omega}(t,k,\eta).
\end{align*}
We saw above that the vorticity experiences enhanced dissipation; since the full formula is a bit long, we will be content with the approximation  $\widehat{\Omega_{\neq}}(t,k,\eta) \sim e^{-c \nu t^3} \widehat{\Omega_{in}}(k,\eta)$. Therefore, the above leads to
\begin{align*}
\mbox{for $k\neq0$}, \quad \widehat{U}(t,k,\eta) \sim -i\begin{pmatrix}kt-\eta \\ k \end{pmatrix} \frac{1}{k^2 + (\eta-kt)^2} \widehat{\Omega_{in}}(t,k,\eta) e^{-c \nu t^3}.
\end{align*}
This formula captures two important effects observed by Orr, one stabilizing, and one destabilizing, together known as the \emph{Orr mechanism}~\cite{Orr07}:
\begin{itemize}
\item[(a)] As $t \to \infty$, it is clear from the above formula that $\widehat{U}_{\neq}$, hence $U_{\neq}$, converges to 0 \textit{uniformly in $\nu$}. 
\item[(b)] On the other hand, the denominator is minimal for $t = \frac{\eta}{k}$, which corresponds to a transient amplification if $\nu t^3 << 1$. 
\end{itemize}
Let us first discuss the damping (a). 
One can be more precise about the rate of decay and obtain the following estimates pointwise-in-time: 
\begin{subequations} \label{ineq:IDlin2D}
\begin{align} 
\| U^1_{\neq}(t) \|_{H^s} & \lesssim \frac{1}{\brak{t}}\norm{\Omega_{in}}_{H^{s+1}} \\
\| U^2_{\neq}(t) \|_{H^s} & \lesssim \frac{1}{\brak{t}^2} \norm{\Omega_{in}}_{H^{s+2}}.
\end{align}
\end{subequations}
Notice crucially that this damping is \emph{independent of Reynolds number}, and indeed, clearly holds also if $\nu = 0$. 
After undoing the coordinate transform \eqref{2dcoords}, this implies similar decay in $L^2$, however, higher $H^s$ norms in general experience slower decay (for $s < 1$) or no decay at all ($s \geq 1$). 
Like enhanced dissipation, the inviscid damping is weak at low wave numbers in $x$. 
In 2D fluid mechanics, the inviscid damping can be easily interpreted as being due to the vorticity being mixed rapidly by the Couette flow, which sends enstrophy to high frequencies. 
Since the Biot-Savart law damps high frequencies, these become less relevant to the velocity field, and eventually, only the shear flow remains.  
The mixing of the vorticity is often called `filamentation' and we remark that it results in a kind of inverse cascade of energy, as the energy in the $x$ dependent modes disappears (in a finite energy, nonlinear problem, this energy would necessarily move into the $x$-independent modes; an analogous effect happens in Landau damping which causes plasmas to heat up as they undergo Landau damping. 

\begin{figure}[ht]
        \centering 
        \includegraphics[width= 0.70 \textwidth]{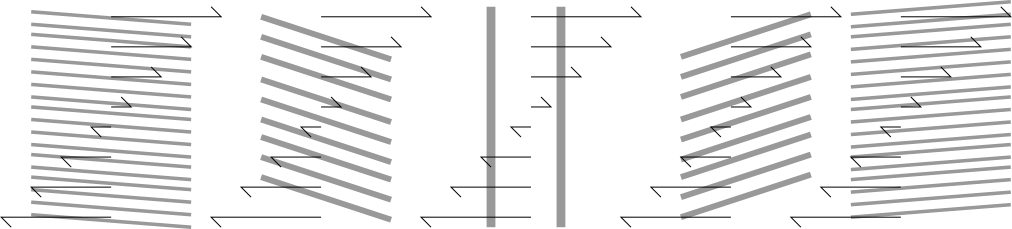}
        \caption{A mode-by-mode visualization of the Orr mechanism: the arrows represent the background flow, while the stripes are the level sets of the function $e^{ik(x-ty) + i\eta y}$ with $\eta/k \gg 1$. Time increases from left to right, and the center image is the critical time $t = \eta/k$. The full linearized solution is simply a super-position of these tilting waves.} \label{fig:Orr}
\end{figure}

In regards to the transient amplification alluded to in (b) above, observe that before the onset of enhanced dissipation, the vorticity $\Omega$ is essentially constant, so that
$$
\mbox{if $k \neq 0$ and $\nu t^3 << 1$,} \quad \frac{\widehat{U}(t=\frac \eta k,k,\eta)}{\widehat{U_{in}}(k,\eta)} \sim \frac{|\eta,k|}{|k|}.
$$
Therefore, if $|\eta| >> |k|$, the velocity field is amplified by a large factor between $t=0$ and the critical time $t=\frac{\eta}{k}$.
In physical terms, this transient growth is due to the fact that the mode of the vorticity in question is initially well-mixed, and then proceeds to \emph{unmix} under the Couette flow evolution. See figure \eqref{fig:Orr} for how this mixing/un-mixing effect appears on each Fourier mode of the vorticity. 
The relevance of the Orr mechanism to hydrodynamic stability has been debated over the years; see e.g. \cite{Orr07,Boyd83,Lindzen88} and \cite{Yaglom12} for a detailed account of how the literature on the topic developed. The importance of the Orr mechanism for nonlinear problems is discussed further in \S\ref{TM}. 

\section{Nonlinear dynamics in dimension $d = 2$}

\subsection{Shear flows as metastable states} \label{sec:shear2Dstab}
%Shear flows are stationary solution of the Euler equation~\eqref{E} of the form 
%\begin{equation}
%\label{shearflow}
%v(y) = \displaystyle \begin{pmatrix} f(y) \\ 0 \end{pmatrix},
%\end{equation}
%for a function $f$.

As we just saw on  the linear level, nonzero frequencies in $x$ ($k \neq 0$) of the velocity field are suppressed (in $L^2$ topology) via: 
\begin{itemize}
\item[(a)] inviscid damping, with a polynomial decay on a time scale $\sim 1$ (with respect to $\nu$);
\item[(b)] enhanced dissipation, with an exponential decay on a time scale $\sim \nu^{-1/3}$.  
\end{itemize}
This implies that, at the linear level, the dependence of $u$ on $x$ is rapidly erased in $L^2$; by the divergence free condition, this implies that $u$ approaches a shear flow rapidly in $L^2$. Convergence of $u$ to a shear flow in higher norms follows after the enhanced dissipation time-scale $\nu^{-1/3}$. 
 However, damping of modes $k=0$ is felt on the much longer time scale $\sim \frac{1}{\nu}$ (which ensures at least the ultimate convergence of the perturbation $u$ to $0$). 
As a conclusion, on a linear level, shear flows are `metastable states', or `intermediate attractors', for the velocity field in $L^2$ on the time range
$$
1 \ll t \ll \nu^{-1}. 
$$
More specifically, for $u_{MS}(y) = (u_0^1(y),0)$, there holds
\begin{align*}
\norm{u(t) - u_{MS}}_{L^2} \lesssim \frac{\epsilon}{\brak{t}}, \quad\quad \textup{for} \;\; 1 \ll t \ll \nu^{-1}. 
\end{align*}
Moreover, we have the following estimate for all $s > 0$ for some $c > 0$,
\begin{align*}
\norm{u(t) - u_{MS}}_{H^s} \lesssim \epsilon \brak{t}^{s-1} e^{-c\nu t^3}, \quad\quad \textup{for} \;\;  1 \ll t \ll \nu^{-1}. 
\end{align*}
Both of these estimates are sharp for the linear problem. Notice that for large $s$ there is a large transient growth before the decay ultimately dominates. 
This is a result of the transfer of enstrophy from low to high frequencies due to the mixing. 

What about the nonlinear problem? Shear flows solve the 1D heat equation for $\nu > 0$ and are stationary solutions for $\nu = 0$, but are they metastable states?
A positive answer to this question is one of the main results in \cite{BMV14,BVW16} (and in some sense \cite{BM13}).  
One of the most important differences between the linear and nonlinear dynamics is the fact that the perturbation itself causes a small adjustment to the background shear flow. 
Specifically, the zero frequency of the first velocity component $u^1_0(t,y)$ is unaffected  by the enhanced dissipation or inviscid damping and hence decays only on the very long $O(\nu^{-1})$ time-scale. 
It follows that over long time-scales, this adjustment to the shear flow can have a very large effect on the solution. 
As a result, \eqref{2dcoords} is no longer necessarily the natural coordinate transformation (this problem is vaguely analogous to a quasilinear scattering problem in dispersive equations, especially for the 2D Euler equations \cite{BM13}). 
As discussed further below, this is one of several major difficulties in the proofs of \cite{BM13,BMV14}. The situation in 3D is even more complicated.

\subsection{Nonlinear resonances} \label{TM} \textit{The frequency cascade scenario.}
Now we arrive at the most fundamental difficulty in understanding the nonlinear dynamics near shear flows at high Reynolds number. 
Specifically, the nonlinear ``resonances'' associated with the non-normal transient growth in the linear problem.
These are not ``true'' resonances as they are not associated with the spectrum of the linear problem, but rather, with the ``pseudo-spectrum'' (see e.g. \cite{TTRD93,Trefethen2005}).  
We will not dwell on such abstract questions, but rather take a direct, hands-on approach to understanding them. 
These nonlinear interactions are much simpler and much better understood in the 2D case -- in 2D there is essentially one leading order nonlinear resonance, whereas the 3D case is far more complicated and it is not clear precisely which interactions will dominate in which situations.

\begin{figure}[ht]
        \begin{center}$ 
\begin{array}{ll}
        \includegraphics[width= 0.50 \textwidth]{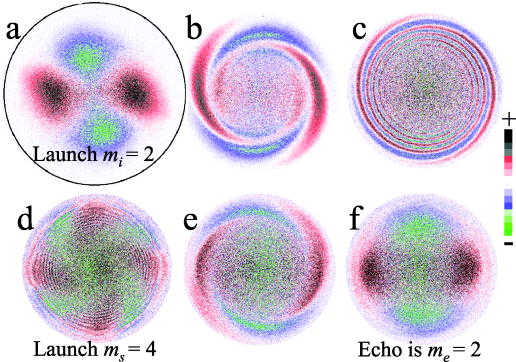} 
\end{array}
$\end{center}
\caption{An experiment showing a hydrodynamic echo in the 2D Euler equations near a radially symmetric vortex. This figure originally appears in \cite{YuDriscoll02}.} \label{figecho}  
\end{figure}

We shall use the coordinate transform conventions defined in \S\ref{PED2}. As we saw in the previous section, the $(k,\eta)$ mode of $\Omega$ (the vorticity in the new variables) has a large effect on $\Phi$ (the streamfunction in the new variables) at the critical time $t = \eta/k$, as can be seen from 
\begin{align*}
\widehat{\Phi}(t,k,\eta) = -\frac{\widehat{\Omega}(t,k,\eta)}{k^2 + \abs{\eta-kt}^2}, 
\end{align*}
(which is simply the Biot-Savart law). 
Physically, if $k\eta >0$, then the mode in the vorticity undergoes a transient un-mixing at the critical time $\eta/k$ before ultimately mixing. 
One then intuitively separates the modes into those for which $\eta/k < t$ and those for which $\eta/k> t$: the un-mixing and mixing modes\footnote{In plasma physics these are sometimes referred to as `phase mixing' and `anti-phase mixing' modes.}.
Then, one imagines a three-wave nonlinear interaction which continues to transfer information from mixing to un-mixing modes to sustain a non-trivial velocity field for long times. 
This effect was first demonstrated in 1968 in collisionless plasmas by the famous \emph{plasma echo} experiments \cite{MalmbergWharton68}. 
Similar `hydrodynamic echoes' were observed in 2D Euler near a radially symmetric vortex via experiments on a pure electron plasma in a strong magnetic field \cite{YuDriscoll02}; 
Figure~\ref{figecho} below is reproduced from that paper.

An \emph{echo} is a nonlinear oscillation in which the velocity field (or electric field in the case of plasmas) decays after an initial disturbance but spontaneously becomes large again later as high-frequency information generated by nonlinear interactions unmixes.
Notice that this is essentially a strong interaction between the linear transient growth coming from the Orr mechanism and the nonlinearity. 
%Such an exchange between the linear transient behavior and nonlinearity is also at the heart of many ideas connected to 3D stability questions; see e.g. \cite{TTRD93,Waleffe95,Chapman02,Trefethen2005} and \S\ref{sec:NulllNonlinear} below. 
Specifically, the scenario we are most concerned with is the following cascade: 
\begin{itemize}
\item Around the critical time $t=\frac{\eta}{k}$, the mode $\widehat{\Phi}(k,\eta)$ is linearly amplified.
\item Through a nonlinear interaction, it causes the mode $\Omega(k-1,\eta)$ to grow...
\item ... which, in turn, linearly amplifies $\widehat{\Phi}(k-1,\eta)$ around time $t = \frac{\eta}{k-1}$.
\item And so on!
\end{itemize}
Solutions with arbitrarily many of such nonlinear plasma echoes have recently been constructed for Vlasov-Poisson \cite{Bedrossian16}, however, the resonances in Vlasov-Poisson are simpler than 2D Euler.  
In order to understand and estimate more quantitatively the growth which can result from this scenario in 2D Euler, we will next develop a toy model.
Closely related studies have been done previously for the 2D Euler equations \cite{VMW98,Vanneste02}, and in some ways our toy model is simpler, but also better suited to getting mathematically rigorous estimates on the dynamics.  

\bigskip
\noindent \textit{The toy model}
The echoes involve the interaction of a mode which is well-mixed and a mode which is near the critical time. In terms of $\Omega$, the vorticity in the shifted coordinates,
 the well-mixed mode will be at a low wave-number and the critical mode will be at a high wave-number (since we are interested in long times). 
Hence, let us formally consider solutions to \eqref{pinguin} of the form $\Omega = \Omega_L + \Omega_H$, where $\Omega_L$ consists only of low frequencies and $\Omega_H$ is a much higher frequency perturbation (but has much less enstrophy).  Let us assume here that neither have contributions to the zero mode in $X$; let us further think of the low frequencies as a background that we are not so interested in, and linearize the equation on high frequencies to obtain
\begin{align*}
\partial_t \Omega_H + \grad^\perp \Phi_{H} \cdot \grad \Omega_L + \grad^\perp \Phi_{L} \cdot \grad \Omega_H  & = \nu \Delta_L \Omega_H \\ 
\Delta_L \Phi_H & = \Omega_H. 
\end{align*}
We argue now that most terms can be dropped in the above: first, we set $\nu = 0$, since the frequency cascade behavior which preoccupies us would  typically occur before dissipation kicks in. Second, since $\nabla^\perp \Phi_L$ is decaying rapidly, and long-term dynamics are of interest, we drop the term $\grad^\perp \Phi_{L} \cdot \grad \Omega_H$. Finally, since $\widehat{\Phi}_H(t,k,\eta) =-\frac{\widehat{\Omega}_H(t,k,\eta)}{k^2 + \abs{\eta-kt}^2}$, we see that the term involving $\partial_Y \Phi_H$ is the most problematic in $\grad^\perp \Phi_{H} \cdot \grad \Omega_L$, as near the critical time $t \sim \eta/k$, $\widehat{\partial_Y\Phi}_H \sim -kt\widehat{\Omega}_H$.
Hence, we reduce to the (linear in $\Omega_H$) model equation 
\begin{subequations} \label{OmToy2}
\begin{align}
\partial_t \Omega_H - \partial_Y \Phi_H \partial_X \Omega_L & = 0, \\ 
\Delta_L \Phi_H & = \Omega_H. 
\end{align}
\end{subequations}
Taking the Fourier transform of \eqref{OmToy2}, the model becomes 
\begin{align*}
\partial_t \widehat{\Omega}_H(t,k,\eta) & = - \sum_{k' \neq 0} \int \frac{\eta'(k-k')}{(k')^2 + (\eta'-k't)^2} \widehat{\Omega}_{H}(t,k',\eta') \widehat{\Omega}_{L}(k-k',\eta-\eta') \,d\eta'. 
\end{align*}
As $\Omega_L$ is concentrated at low frequencies, let us make the simplifying approximation that $\eta \approx \eta'$ and $k' = k \pm 1$. 
This gives us the following model
\begin{equation}
\partial_t \widehat{\Omega}_H(t,k,\eta) = -\sum_{k' = k \pm 1} \frac{\eta (k-k')}{(k')^2 + (\eta-k't)^2} \widehat{\Omega}_H(t,k',\eta) \Omega_L(k-k',0). \label{ODETM} 
\end{equation}
The idea is that this model might be useful in estimating the leading order, long-time dynamics of how enstrophy in the high frequency perturbation $\Omega_H$ is transferred between modes in $X$ as a result of the interactions with the low frequencies. 
For a general solution, the simplifying assumptions we made in deriving \eqref{ODETM} are all inaccurate, and we will not be able to say that \eqref{ODETM} is necessarily an accurate model for the dynamics (although the echo construction of \cite{Bedrossian16} shows that an analogous toy model for echoes in the Vlasov-Poisson equations can be an accurate approximation; see \S\ref{sec:LandauMixing}). 
However, in \cite{BM13,BMV14}, we instead use \eqref{ODETM} to estimate the \emph{worst case scenario}, that is, we use \eqref{ODETM} to get an upper bound on how much enstrophy can be transferred from low frequencies to high frequencies. See \S\ref{sec:FourerMult} for a discussion of how this is done to enable the proofs in \cite{BM13,BMV14}. 
\section{Linear dynamics in dimension $d=3$}  
\label{PED3}

\subsection{Coordinate transform}
\label{NC3D}
In three dimensions, the vorticity formulation seems less advantageous (especially for the nonlinear problem; see \cite{BGM15I} for a brief discussion).
One reason for this is that, unlike the 2D case in \eqref{def:vortform}, the linear evolution of the vorticity is not a passive transport equation, but instead also involves a non-local term giving rise to vortex stretching.  
Even in the limit $t \rightarrow \infty$, it is not true that the dominant dynamics are passive scalar transport uniformly in $\nu$. 
Nevertheless, it will still be advantageous to begin by factoring out by the transport as in \eqref{2dcoords}. 
Hence, we set
\begin{equation*}
\left\{
\begin{array}{l}
X = x-ty \\
Y = y \\
Z= z.
\end{array}
\right.
\end{equation*}
As above, we capitalize the physical variables taken in this new set of coordinates: $U(X,Y,Z) = u(x,y,z)$, and denote $\nabla_L = \begin{pmatrix} \partial_X \\ \partial_Y - t\partial_X \\ \partial_Z \end{pmatrix}$ and $\Delta_L = \nabla_L \cdot \grad_L$ for the differential operators in the new coordinates. 
In these new coordinates, the equation satisfied by $U$ reads
\begin{equation}
\label{eqUXYZ}
\left\{
\begin{array}{l}
\partial_t U + U \cdot \nabla_L U = \begin{pmatrix} - U^2 \\ 0 \\ 0 \end{pmatrix} - \nabla_L P + \nu \Delta_L U \\
\nabla_L \cdot U = 0 \\
U(t=0) = U_{in}.
\end{array}
\right.
\end{equation}
Linearizing this for small perturbations gives the system  
\begin{equation}
\label{eqUXYZlin}
\left\{
\begin{array}{l}
\partial_t U = \begin{pmatrix} - U^2 \\ 0 \\ 0 \end{pmatrix} - \nabla_L P + \nu \Delta_L U \\
\nabla_L \cdot U = 0 \\
U(t=0) = U_{in}.
\end{array}
\right.
\end{equation}
Notice that in \eqref{eqUXYZlin}, the pressure is given by
\begin{align}
\Delta_L P = -2 \partial_{x}U^2. \label{P3Dlin}
\end{align}
%As the non-trivial terms in \eqref{eqUXYZlin} depend only on $U^2$, it is natural to first deduce the behavior of $U^2$ and then use this information to deduce the behavior of the other two variables. 

\subsection{Inviscid damping and enhanced dissipation of $U^2$} 
\label{IDOM3D}

As it happens, inviscid damping is still relevant in 3D, despite the vortex stretching. 
Interestingly, \cite{BGM15I} is the first work known to the authors that explicitly points out the relevance of inviscid damping to 3D fluid mechanics.  
This is likely because it, at first, appears to be a relatively minor detail in the overall physical picture of the dynamics. 
However, as discussed below in Remark \ref{rmk:ID3D}, the presence of inviscid damping will be absolutely critical to the proof of the results in \cite{BGM15I,BGM15II,BGM15III}. 
In order to see the effect in 3D, we introduce a new variable, $Q^2 = \Delta_L U^2$. 
This unknown is well-known in the hydrodynamic stability of shear flows, see e.g. \cite{SchmidHenningson2001}, and was apparently introduced by Kelvin \cite{Kelvin87} (more precisely, the unknown $q^2 = \Delta u^2$ is normally the unknown studied). 
A small computation reveals that $Q^2$ solves
\begin{align}
\partial_t Q^2 = \nu \Delta_L Q^2. \label{linQ2}
\end{align} 
(this calculation is actually easier to carry out in $(x,y,z)$ and then transfer to $(X,Y,Z)$ afterwards). Since
$$
U^2(t) = \Delta_L^{-1} Q^2(t),
$$
we derive the following, \emph{independent of Reynolds number}, analogous to \eqref{ineq:IDlin2D}, 
\begin{align}
\norm{U^2_{\neq}(t)}_{H^s} \lesssim \frac{1}{\brak{t}^2}\norm{Q^2_{in}}_{H^{s+2}}. \label{ineq:3DID}
\end{align}
Similarly, the Orr mechanism manifests on $U^2$ similar to the way it manifests in 2D. 
By an analogous computation using \eqref{P3Dlin}, we also have the rapid decay of the pressure:   
\begin{align}
\norm{P(t)}_{H^s} & \lesssim \frac{1}{\brak{t}^4}\norm{\Omega_{in}}_{H^{s+4}}.  \label{ineq:3DIDpres}
\end{align}
Finally, since \eqref{linQ2} is the same as \eqref{linOmega2D}, it is clear that $Q^2$ experiences the same enhanced dissipation effect as $\Omega$ does in 2D (see \S\ref{sec:ED2D}).

\subsection{Vortex stretching and enhanced dissipation of $U^{1,3}$}
>From \eqref{ineq:3DID}, \eqref{ineq:3DIDpres}, we see that the $P$ and $U^2$ terms in 
\begin{align*}
\partial_t \begin{pmatrix} U^1 \\ U^3 \end{pmatrix} & = \begin{pmatrix} -U^2 \\ 0 \end{pmatrix} - \begin{pmatrix} \partial_X P \\ \partial_Z P \end{pmatrix} + \begin{pmatrix} \nu \Delta_L U^1 \\ \nu \Delta_LU^3 \end{pmatrix} \\ 
0 & = \partial_X U^1 + (\partial_Y - t\partial_X)U^2 + \partial_Z U^3, 
\end{align*}
decay rapidly on $O(1)$ time-scales ($O(1)$ relative to the Reynolds number). It is therefore acceptable to choose $u^2_{in} = 0$ (this does not alter the generic behavior), for which on time-scales $1 \lesssim t$, the above becomes 
\begin{align*}
\partial_t \begin{pmatrix} U^1 \\ U^3 \end{pmatrix} & =  \begin{pmatrix} \nu \Delta_L U^1 \\ \nu \Delta_LU^3 \end{pmatrix} \\ 
\partial_X U^1 +\partial_Z U^3 & = 0. 
\end{align*}
There are two main observations here: (A) $U^1_{\neq}$ and $U^3_{\neq}$ experience the same enhanced dissipation as $U^2$ but (B) they do \emph{not} experience any inviscid damping. 
In fact, for time-scales $1 \ll t \ll \nu^{-1/3}$, the definition of the coordinate transform gives that for some fixed, \emph{time-independent} $U^{1,\infty}_{\infty}$
\begin{align*}
u^{1}(t,x,y,z) & \sim U^1_\infty(x-ty,y,z) \\ 
u^{3}(t,x,y,z) & \sim U^3_\infty(x-ty,y,z). 
\end{align*}
See \cite{BGM15I}. Hence, in general, the vorticity $\omega = \grad \times u$ grows (e.g. in $L^2$) linearly in time as $O(\epsilon \brak{t})$ until the dissipation time-scale $\nu^{-1/3}$ and $u^{1,3}$ experience a matching linear-in-time \emph{kinetic energy} cascade. 
This vortex stretching adds extra complexity to both the proofs and the actual nonlinear dynamics in 3D vs 2D (for example, see the complexity difference between the 2D work \cite{BVW16} and the analogous 3D work \cite{BGM15III}). 
Moreover, the Orr mechanism in 3D not only involves the un-mixing of enstrophy from small-to-large scales as originally observed by Orr, but now also involves a more subtle vorticity stretching effect which is anisotropic in $(X,Y)$ frequencies (see e.g. \cite{PradeepHussain06}). 
This will ensure that the 3D resonances are significantly more complicated than 2D.

\subsection{The lift up effect} 
There is an additional linear effect, which was discovered in \cite{EllingsenPalm75} (see also \cite{landahl80} for extensions to low wave-numbers in $x$). 
This effect is active for $x$-independent modes ($U_0$): for these modes, $\nabla_L = \nabla$ and $\Delta_L = \Delta$, so that the linearization is given by 
\begin{equation*}
\partial_t U_0  = \begin{pmatrix}- U^2_0 \\ 0 \\ 0 \end{pmatrix} + \nu \Delta U_0. 
\end{equation*}
Fortunately, this system can be solved explicitly:
\begin{equation*}
\left\{
\begin{array}{l}
U^1_0(t) = e^{\nu t \Delta} \left[ U^1_{in,0} - t U^2_{in,0} \right] \\
U^2_0(t) = e^{\nu t \Delta} U^2_{in,0} \\
U^3_0(t) = e^{\nu t \Delta} U^3_{in,0}.
\end{array}
\right.
\end{equation*} 
This means that $u^1$ undergoes a linear growth until time $\sim \frac{1}{\nu}$, at which point the viscosity will begin to slowly relax $u^1$ to zero.  
In fact, the dynamics are essentially the same as the canonical non-normal ODE system: 
\begin{align*}
\partial_{t} \begin{pmatrix} X_1 \\ X_2 \end{pmatrix} = \begin{pmatrix} -\nu X_1 - X_2 \\ -\nu X_2 \end{pmatrix}. 
\end{align*}
This transient growth is called the \emph{lift-up effect} and is one of the main destabilizing mechanisms in three dimensions (it is absent in dimension two since there the divergence free condition imposes $U^2_0 = 0$). 
This instability is one of the primary culprits for subcritical transition in 3D and implies that 3D flows generally have very different dynamics from their 2D counter-parts. For example, one can easily verify that \emph{every} non-trivial 3D shear flow undergoes this transient growth (this is easily seen on the 3D Euler equations, $\nu = 0$, where the lift-up effect becomes an unbounded, linear-in-time algebraic instability). 
The lift-up effect and the nonlinear effects of vorticity stretching together imply that, at high Reynolds numbers, the hydrodynamic stability of 2D shear flows is very different from that of 3D shear flows (consistent with experimental observations). 

\section{Nonlinear dynamics in dimension three}

\subsection{Streaks} \label{sec:streaks}
Observe that for $u_E(y) = (y,0,0)^t$, if the initial data in \eqref{eq:perteqns} is independent of $x$, then this remains true for all time: $u(t,x,y,z) = u(t,y,z)$. 
In fact, $u^{2,3}(t)$ solves the 2D Navier-Stokes equations in $(y,z) \in \Real \times \Torus$ whereas
$u^1$ solves the following linear advection-diffusion equation (solutions to Navier-Stokes of this type are generally called `2.5 dimensional'; see e.g. \cite{MajdaBertozzi}):  
$$
\partial_t u^1  + (u^2 \partial_y + u^3 \partial_z) u^1 = -u^2 +  \nu \Delta u^1. 
$$
We will refer to these solutions as \emph{streaks}.
Due to the lift-up effect, the vast majority of the kinetic energy of a streak is in $u^1$ \cite{ReddySchmidEtAl98}; in experiments and computer simulations this, together with the low wave-number in $x$, gives them a very distinct streaky appearance \cite{klebanoff1962,bottin98,ReddySchmidEtAl98,SchmidHenningson2001} (for $x \in \Real$, the energy is piling into smaller and smaller wave-numbers in $x$; see \cite{landahl80}). 
Using the lift-up effect, it is easy to construct streaks which are initially $O(\nu)$ in $L^2$ but at $t \approx \nu^{-1}$, satisfy $\norm{u^1(t)}_{L^2} \gtrsim 1$.
As solutions to the 2D Navier-Stokes are global in time, this is a class of global solutions to the 3D Navier-Stokes equations which are from equilibrium relative to $\nu$.

If one only considers the linear problem, then recall that $U_{\neq}$ is dissipated on the short time scale $\nu^{-1/3}$ so that streaks are potentially meta-stable states of the 3D nonlinear problem somewhat analogous to the shear flows in 2D. 
There are some important differences however, for example, due to the vortex stretching, the velocity field is \emph{far} from the streak solutions in $L^2$ for $t \ll \nu^{-1/3}$, whereas in 2D, the inviscid damping ensures that the velocity is close to a shear flow after times which are $O(1)$ relative to $\nu$ (a difference with major nonlinear implications). 
Assuming that the solution does not transition to a turbulent state, for $t \gg \frac{1}{\nu}$, dissipation in the $k=0$ modes send $U_0$ to zero. 

\subsection{A scenario for subcritical transition}
If streaks are meta-stable states near the stability threshold, a possible storyline emerges: 
%The linear, as well as nonlinear, heuristics, lead to the following storyline, which we are able to confirm analytically in Gevrey topology (the picture is less clear in Sobolev topology):

\begin{itemize}
\item Start at $t=0$ with a perturbation $u_{in}$ of size $\sim \epsilon$.
\item At times $\displaystyle t \gg \frac{1}{\nu^{1/3}}$, enhanced dissipation kicks in and dissipates $U_{\neq}$. As a result, the solution looks essentially like a streak.
\item At time $\displaystyle t \sim \frac{1}{\nu}$, the maximal linear growth from the lift up effect predicts $U_0^1$ could have size $\displaystyle \sim \frac{\epsilon}{\nu}$.
\item There are now two possibilities
\begin{itemize}
\item If $\displaystyle \frac{\epsilon}{\nu} \ll 1$, the maximal size of the streak is $\frac{\epsilon}{\nu}$ and the dissipation of modes $k=0$ should, for times $\displaystyle t \gg \frac{1}{\nu}$, bring $u$ back to zero without transitioning to a fully nonlinear state.  
\item If $\displaystyle \frac{\epsilon}{\nu} \gg 1$, we eventually leave the perturbative regime, and additional instabilities should be expected at around the time $\epsilon^{-1}$, when the streak becomes $O(1)$. 
\end{itemize}
\end{itemize}
If this storyline is accurate, it would lead to two conclusions: we could determine the stability threshold is precisely $\gamma = 1$ and we would have that all instabilities would result from the `secondary instability' of a streak. 
That is, a streak grows such that $\norm{u_0^1}_{L^2} = O(1)$, and then the slowly varying shear flow $(y + u_0^1(t),0,0)$ develops a true exponential instability; e.g. if one formally freezes $t$ \cite{ReddySchmidEtAl98} then one can (formally) create an inviscid, inflection-point shear flow instability. Such an instability would develop very rapidly relative to the time-scales on which the streak varies, and carry the solution rapidly into a strongly nonlinear regime (though it is not quite accurate to suggest it will necessarily go straight into a turbulent state; see \cite{DuguetEtAl10} and the references therein for more details on this). 
This secondary instability is called \emph{streak breakdown} in the fluid mechanics literature, and has attracted a lot of attention as it  %(see \cite{bottin98,ReddySchmidEtAl98,Chapman02,DuguetEtAl2010,SchmidHenningson2001} and the references therein). 
is the scenario most often observed in computer simulations and physical experiments (see \cite{TTRD93,bottin98,ReddySchmidEtAl98,Chapman02,SchmidHenningson2001,DuguetEtAl2010} and the references therein).

The works \cite{BGM15I,BGM15II} together confirm at least the beginning of this general picture in Gevrey-$s^{-1}$ with $s > 1/2$, that is  $\norm{e^{\lambda \brak{\grad}^s}u_{in}}_{2} = \epsilon$ and for $\epsilon \lesssim \nu^{2/3+\delta}$ for any small $\delta > 0$.
Note this latter condition is simply a requirement that the data is not `too far' above the stability threshold; it is obvious we need an assumption of this general form, however it is an open question as to whether $2/3$ should be the sharp exponent. 
The two conditions are used to deduce that fully 3D nonlinear effects are dominated by the enhanced dissipation.  

\subsection{Null forms and nonlinear interactions} \label{sec:NulllNonlinear}
The scenario that was just proposed is relevant if the --mostly linear-- heuristics which were derived in the previous paragraphs apply. 
This can only be the case if nonlinear effects remain small; one reason why this is the case in the regimes we consider goes by the name of \textit{null form}. The idea is the following: the most potentially threatening nonlinear interactions between linear waves are ruled out by the particular structure of the equation.
These `weakly nonlinear' interactions are significantly more complicated in 3D than the 2D echoes discussed in \S\ref{TM}, however, the nonlinear structure of the 3D Navier-Stokes equations near Couette is still very specific and we can recognize a few special structures. 
For example, here are a few basic observations (not at all exhaustive): 
\begin{itemize}
\item[(a)] Self-interactions of the lift up: recall that $U_0^1$ grows linearly until time $\sim \frac{1}{\nu}$. A quadratic interaction of $U_0^1$ with itself would be catastrophic, but it is not allowed by the nonlinearity of the equation.
\item[(b)] Inviscid damping vs convection: Recall that the convection term above in \eqref{eqUXYZ} reads, when spelled out
$$
U \cdot \nabla_L U = [U^1 \partial_X + U^2 (\partial_Y - t\partial_X) + U^3 \partial_Z] U.
$$
The second summand, $U^2 (\partial_Y - t\partial_X) U$, seems quite threatening, since it contains a linearly growing factor. Fortunately, it is paired with $U^2$ which, as we saw in Section~\ref{IDOM3D}, is damped independently of Reynolds number via inviscid damping. A similar observation also damps problematic terms in the pressure as a similar $U^2 (\partial_Y - t\partial_X) U^j$ structure is preserved therein as well. 
\end{itemize}
\begin{remark} \label{rmk:ID3D}
Note that observation (b) emphasizes the key relevance of inviscid damping for understanding the nonlinear problem in 3D. Indeed, inviscid damping might look like a minor detail in the statements of our results \cite{BGM15I,BGM15II,BGM15III} or in the general physical picture of the dynamics laid out above. 
However, inviscid damping plays a \emph{crucial} role in the \emph{proof} due to this kind of nonlinear structure. 
\end{remark} 

There are several nonlinear mechanisms which have the potential to cause instability and many have been proposed as important in the applied mathematics and physics literature for understanding transition, see e.g. \cite{Craik1971,TTRD93,ReddySchmidEtAl98,SchmidHenningson2001} and the references therein. 
We are particularly worried about so called ``bootstrap'' mechanisms \cite{TTRD93,VMW98,Vanneste02,Trefethen2005,Waleffe95,BM13}: nonlinear interactions 
that repeatedly excite growing linear modes. %, a kind of nonlinear ``pseudo-resonance''. 
We will classify the main effects by the $x$ frequency of the interacting functions: denote for instance $0 \cdot \neq \,\to\, \neq$ for the interaction of a zero mode (in $x$) and a non-zero mode (in $x$) giving a non-zero mode (in $x$), and similarly, with obvious notations, $0 \cdot 0 \to 0$, $\neq \cdot \neq \, \to \,\neq$, and $\neq \cdot \neq \,\to 0$.

\begin{itemize} 
\item ($0 \cdot 0 \to 0$)  These correspond to self-interactions within the streak. 
\item ($0 \cdot \neq \,\to\, \neq$) These interactions are essentially those that arise when linearizing an $x$-dependent perturbation of a streak and so are connected to the secondary instabilities observed in larger streaks \cite{ReddySchmidEtAl98,Chapman02}. The instabilities which are most commonly observed in experiments are generally related to secondary linear instabilities of the streaks, and so it is unsurprising that most of the leading order interactions are of this type. 
%or \emph{secondary instability}, this effect is the transfer of momentum from the large $u_0^1$ mode to other modes (actually $u_0^2$ and $u_0^3$ also matter, but less). 
%Due to the special `skew' structure of the nonlinearity, as seen in the 2.5 dimensional nature of the streaks in \eqref{def:streak}, this effect can only occur on modes which are $x$-dependent. 
%hese involve a zero frequency and non-zero frequency $k$ interacting to amplify the same mode $k$, or the $k$ mode of a different component, e.g. $u_0^1$ and $u_k^3$ together force $u_k^2$. 
\item ($\neq \cdot \neq \,\to \,\neq$) These effects include the 3D variants of the 2D hydrodynamic echo phenomenon as observed in \cite{YuDriscoll02,YuDriscollONeil}: nonlinear interactions of $x$-dependent modes forcing unmixing modes \cite{Morrison98,Vanneste02,BM13} -- a nonlinear manifestation of the Orr mechanism. 
In 3D, this is still important and the range of possible interactions  is much wider (see e.g. \cite{Craik1971,SchmidHenningson2001,Yaglom12}). 
%is connected not only to unmixing but also to the vortex stretching. 
%This involves two non-zero frequencies $k_1$, $k_2$ interacting to force mode $k_1 + k_2$ with $k_{1},k_2,k_1 + k_2 \neq 0$. %As above, these interactions can of course involve modes from different components of the solution. %Since these interactions are ultimately due to a fundamentally 2D unmixing effect, we can consider these terms as ``2.5 dimensional echoes''. 
%Since these involve the interaction of only non-zero frequencies, they should only be problematic for short times: for $t \gtrsim \nu^{-1/3}$, this effect should be wiped out by the enhanced dissipation. 
\item ($\neq \cdot \neq \,\to 0$) These effects are the nonlinear feedback from $x$-dependent modes back into $x$-independent modes. %This can occur directly into $u_0^1$ or it can occur in $u_0^2$, which due to the lift-up effect, has a strong effect on $u_0^1$. Forcing into $u_0^3$ can also have an effect as this can subsequently amplify $u_0^2$, which feeds back into $u_0^1$.
%This involves two non-zero frequencies $k$ and $-k$ interacting to force a zero frequency (in general this could involve a variety of the components). Similar to \textbf{(3DE)}, this effect is over-powered by the enhanced dissipation after $t \gtrsim \nu^{-1/3}$. 
\end{itemize} 
All of these interactions are coupled to one another, and one can imagine bootstrap mechanisms involving several of them (e.g. $u_0^1$ forces a non-zero mode which unmixes and then strongly forces $u_0^2$ which strongly forces $u_0^1$ via the lift-up effect and repeat). 

Toy models such as that derived in \S\ref{TM} above can be written down for 3D, however, it is not really practical to simultaneously include all of the leading order interactions in the same model. 
Instead, in \cite{BGM15I,BGM15II}, we were content with deriving a much rougher toy model meant only to help provide upper bounds on the nonlinear resonances. 
Even then, the toy models are 6x6 ODEs involving many terms which contain information involving all of the linear dynamics deduced above.  

\section{Statement of the theorems} \label{sec:Statements}

%We state below the theorems obtained in Gevrey regularity, which are optimal as far as the asymptotic stability exponent (which we called $\gamma$) goes:
%in dimension $2$, $\gamma = 0$, and indeed, the Couette flow is stable in Euler; while in dimension $3$, $\gamma =1$ is optimal, see below.

For the sake of simplicity, we focus in this section on results in Gevrey regularity, and only give short statements of a selection of theorems, referring the interested reader to the original papers.

\subsection{The case of dimension 2}

\begin{thm} [Stability of Couette in dimension 2~\cite{BM13}~\cite{BMV14}]
Consider~\eqref{NS} for $\nu>0$, or~\eqref{E}, in which case $\nu = 0$. Fix $s \in (\frac{1}{2},1)$ and $\lambda_0 > \lambda^\prime > 0$. Then, if
$$
\norm{u_{in}}_{\mathcal{G}^{\lambda;s}} = \epsilon
$$
is sufficiently small (depending only on $s,\lambda_0,\lambda'$), then the unique, classical solution $u(t)$ to \eqref{eq:perteqns} with initial data $u_{in}$ is global in time and enjoys the following estimates (all implicit constants might depend on $s$, $\lambda_0$, $\lambda'$, but not on $\nu$)
\begin{itemize}
\item[(i)] inviscid damping and enhanced dissipation of the velocity field, 
$$
\norm{u^1_{\neq}(t)}_{L^2} + \brak{t}\norm{u^2(t)}_{L^2}  \lesssim \frac{\epsilon}{\brak{t} \brak{\nu t^3}^{10}}.
$$
\item[(ii)] enhanced dissipation on time scales $\gtrsim \frac{1}{\nu^{1/3}}$:
$$
\norm{\omega_{\neq}(t,x + ty + t\psi(t,y))}_{\G^{\lambda^\prime;s}}  \lesssim \frac{\epsilon}{\brak{\nu t^3}^{10}},
$$
(here $\psi(t,y)$ is an $O(\epsilon)$ correction to the mixing which depends on the disturbance).
\item[(iii)] classical viscous decay of the zero modes
$$
\norm{u_{0}}_{L^2} \lesssim \frac{\epsilon}{\langle \nu t \rangle^{1/4}}.
$$
\end{itemize}
\end{thm}
\begin{remark}
One can also make more precise the meta-stability assertions in \S\ref{sec:shear2Dstab} by showing the existence of a shear flow $u_\infty = (u_\infty^1(y),0)$ such that $u(t) \approx u_{\infty}$ on time-scales $1 \ll t \ll \nu^{-1/3}$. Similarly, one can prove the existence of an $\omega_\infty = \omega_\infty(x,y)$ such that $\omega(t,x,y) \approx \omega_\infty(t,x+ty + t\psi(t,y)) $ for times $1 \ll 1 \ll \nu^{-1/3}$, which is analogous to scattering in dispersive equations. See \cite{BM13,BMV14} for more discussions.  
\end{remark}

To summarize, the above theorem states that the set of shear flows is asymptotically stable in Gevrey topology in the Euler equations (in a sense), and, uniformly in $\nu>0$, in the Navier-Stokes equations. 
When viscosity is turned on, an additional effect appears, namely enhanced dissipation, but this is a stabilizing mechanism and hence the result is uniform in $\nu$.

\subsection{The case of dimension 3}

\begin{thm}[Below threshold dynamics \cite{BGM15I}] \label{thm:Threshold} 
Fix $s \in (1/2,1)$, and $\lambda_0 > \lambda^\prime > 0$.
Then there exists $c_0$ such that: if
\begin{align} 
\norm{u_{in}}_{\mathcal{G}^{\lambda;s}} = \epsilon < c_0 \nu, \label{ineq:QuantGev2}
\end{align} 
then the unique, classical solution $u(t)$ to \eqref{eq:perteqns} with initial data $u_{in}$ is global in time and enjoys the following estimates (all implicit constants might depend on $s$, $\lambda_0$, $\lambda'$, but not on $\nu$)
\begin{itemize}
\item[(i)] the rapid convergence to a streak through enhanced dissipation:
\begin{subequations} \label{ineq:uidamping}
\begin{equation} 
%\norm{u_{\neq}(t,x + ty + t\psi(t,y,z),y,z)}_{\G^{\lambda^\prime;s}}  \lesssim \frac{\epsilon \brak{t}^{\frac{1}{10}}}{\brak{\nu t^3}^{10}}. \label{ineq:u1damping}
\norm{u_{\neq}(t)}_{L^2}  \lesssim \frac{\epsilon \brak{t}^{\frac{1}{10}}}{\brak{\nu t^3}^{10}}. \label{ineq:u1damping}
\end{equation}  
\end{subequations}
%(here $\psi(t,y,z)$ is an $O(c_{0})$ correction to the mixing which depends on the disturbance). 
\item[(ii)] transient growth of the streak for $t < \frac{1}{\nu}$  through the lift-up effect: 
\begin{subequations} \label{ineq:trans}
\begin{align}
\norm{u^1_0(t) -  \left(e^{\nu t \Delta} \left(u^{1}_{in \; 0} - t u_{in \; 0}^2\right) \right) }_{\G^{\lambda^\prime;s}} & \lesssim \left( \frac{\epsilon}{\nu} \right)^2 \label{ineq:u01grwth1}\\
\norm{u^{2}_0(t) -  e^{\nu t \Delta} u^{2}_{in \; 0} }_{\G^{\lambda^\prime;s}} + \norm{u^{3}_0(t) -  e^{\nu t \Delta} u^{3}_{in \; 0} }_{\G^{\lambda^\prime;s}} & \lesssim \left( \frac{\epsilon}{\nu} \right)\epsilon \label{ineq:u023} 
\end{align}
\end{subequations} 
\item[(iii)] decay of the background streak for $t > \frac{1}{\nu}$:
\begin{equation} 
\norm{u(t)}_{\G^{\lambda^\prime;s}} \lesssim \frac{\epsilon}{\nu \brak{\nu t}^{1/4}}.
\end{equation}
\end{itemize}
\end{thm} 
 
In particular, \eqref{ineq:u01grwth1} (together with the other estimates) show that it is possible to find initial data of the size $\eps = c_0\nu$ but such that at some finite time $t_\star \approx \nu^{-1}$, $\norm{u_0^1(t)}_{L^2} \gtrsim c_0$. 
Hence, one has arbitrarily small solutions which become $O(1)$ with respect to the Reynolds number. 
However, these solutions are not quite large enough to trigger transition, as is evident by the fact that all the estimates in Theorem \ref{thm:Threshold} continue till $t \rightarrow \infty$. 

\bigskip

Next, one is interested in studying the dynamics of solutions above the stability threshold. 
This is the content of the next theorem. 
 
\begin{thm}[Above threshold dynamics \cite{BGM15II}] \label{thm:SRS}
Fix $s \in (1/2,1)$, and $\lambda_0 > \lambda^\prime > 0$. Then, for all $\delta > 0$ all $\nu$ and $c_0$ sufficiently small (depending only on $s,\lambda_0,\lambda',\delta$), if  
\begin{align} 
\norm{u_{in}}_{\mathcal{G}^{\lambda;s}} = \epsilon < \nu^{2/3-\delta}, \label{ineq:QuantGev2}
\end{align} 
then the unique, classical solution $u(t)$ to \eqref{eq:perteqns} with initial data $u_{in}$ exists at least until time $T_F = c_0 \epsilon^{-1}$ and enjoys %the following estimates (all implicit constants might depend on $s$, $\lambda_0$, $\lambda'$, but not on $\nu$)
%\begin{align} 
%\norm{u_S}_{\mathcal{G}^{\lambda_0;s}} + e^{K_0\nu^{-\frac{3s}{2(1-s)}}}\norm{u_R}_{H^{3}} & \leq \epsilon, \label{ineq:QuantGev22}
%\end{align} 
%then the unique, classical solution to \eqref{eq:perteqns} with initial data $u_{in}$ exists at least until time $T_F = c_0 \epsilon^{-1}$ 
 the following estimates with all implicit constants independent of $\nu$, $\epsilon$, $c_0$ and $t$:
\begin{itemize}
\item[(i)] the rapid convergence to a streak  through enhanced dissipation: 
\begin{equation} 
%\norm{u_{\neq}(t,x + ty + t\psi(t,y,z),y,z)}_{\G^{\lambda^\prime;s}} \lesssim \frac{\epsilon t^{\frac{1}{10}}}{\brak{\nu t^{3}}^{10}} \label{ineq:u2damping2}
\norm{u_{\neq}(t)}_{L^2} \lesssim \frac{\epsilon t^{\frac{1}{10}}}{\brak{\nu t^{3}}^{10}} \label{ineq:u2damping2}
\end{equation} 
\item[(ii)] transient growth of the streak for $t < T_F$ through the lift-up effect:
\begin{align}
\norm{u^1_0(t) -  e^{\nu t\Delta}\left(u_{in \; 0}^1 - tu_{in \; 0}^2\right) }_{\G^{\lambda^\prime;s}} & \lesssim c_0^2 \label{ineq:u01grwth} \\ 
\norm{u^2_0(t) -  e^{\nu t\Delta} u_{in \; 0}^2}_{\G^{\lambda^\prime;s}} + \norm{u^3_0(t) -  e^{\nu t\Delta} u_{in \; 0}^3}_{\G^{\lambda^\prime;s}} & \lesssim c_0 \epsilon; \label{ineq:u023Duhamel}
%\left(e^{\nu t \Delta} u^{1}_{0 \; in} - \int_0^t e^{\nu (t-\tau) \Delta} u_{0 \; in}^2 d\tau \right)
\end{align}
\item[(iii)] uniform control of the background streak for $t < T_F$:
\begin{subequations}  
\begin{align} 
\frac{1}{\brak{t}} \norm{u^1_0(t)}_{\G^{\lambda^\prime;s}}+
\norm{u^2_0(t)}_{\G^{\lambda^\prime;s}} + \norm{u^3_0(t)}_{\G^{\lambda^\prime;s}} \lesssim \epsilon; 
\end{align}
\end{subequations}  
(where $\psi(t,y,z)$ is an $O(\epsilon t)$ correction to the mixing which depends on the disturbance).
\end{itemize}
\end{thm} 
The key observation above is that Theorem \ref{thm:SRS} shows that near the transition threshold, the \emph{only} possible instability is the secondary instability of a streak as envisioned in \cite{TTRD93,ReddySchmidEtAl98,Chapman02} and others. 
%That is, a streak grows such that $\norm{u_0^1}_{L^2} = O(1)$, and then the slowly varying shear flow $(y + u_0^1(t),0,0)$ develops a true exponential instability; e.g. if one formally freezes $t$ \cite{ReddySchmidEtAl98} then one can, at least formally, create an inviscid, inflection-point shear flow instability. Such an instability would develop very rapidly relative to the time-scales on which the streak varies, and carry the solution rapidly into a nonlinear regime (though it is not quite accurate to suggest it will necessarily go straight into a turbulent state, as discussed above). 

\section{Rigorous mathematical proofs} \label{sec:Rigor}

\subsection{Nonlinear change of coordinates} \label{sec:nonChange} 
We saw in Sections~\ref{NC2D} and~\ref{NC3D} that a change of variable was necessary to take into account the effect of convection by the Couette flow in the linearized equation.

\bigskip

\noindent \textit{The 2D case}. In the work \cite{BVW16}, the viscosity is (barely) large enough to imply that the viscous effects dominate the adjustments to the shear flow over long times, and hence in \cite{BVW16}, \eqref{2dcoords} is essentially sufficient. However, in \cite{BM13,BMV14} the viscosity is non-existent or arbitrarily small relative to the size of the perturbation (respectively), and hence, the perturbation to the shear cannot be neglected. 
In these cases, we use the following ansatz with a function $\psi$ to be determined:
\begin{subequations} \label{def:nonlin2Dcoord}
\begin{align}
X & = x - ty - t\psi(t,y) \\ 
Y & = y + \psi(t,y).
\end{align}
\end{subequations}
The purpose of $x \mapsto X$ is to account for the mixing due to the background shear flow and, equivalently, to find a coordinate system in which the vorticity could have uniform-in-time higher regularity estimates. 
The purpose of $y \mapsto Y$ is so that $\partial_y \mapsto (1 + \partial_y \psi) (\partial_Y - t\partial_X)$, which implies that the critical times still occur at $t \sim \eta/k$ (see \eqref{def:SymDeltat} and \cite{BM13}). 
This ansatz is applied and we derive the following: 
\begin{align}
\partial_t \Omega +  \begin{pmatrix} u^1 - \partial_t(t\psi) \\ u^2 + \partial_t \psi - \nu \partial_{yy}\psi \end{pmatrix} \cdot \begin{pmatrix}\partial_X\Omega \\ (1+\partial_y \psi)(\partial_Y - t\partial_X) \Omega \end{pmatrix}  = \nu \partial_{XX}\Omega + \nu(1+\partial_y \psi)^2(\partial_Y - t\partial_X)^2 \Omega. \label{eq:midstepOm}
\end{align} 
Note that $\psi$ and $u$ are still in $(x,y)$ coordinates. 
If one assumes that the linear problem is a good approximation for the dynamics, then $u_0^1$ is the only contribution to the drift which is not decaying quickly. 
Hence, it makes sense to choose $\psi$ in order to eliminate this, which is the choice made in \cite{BM13,BMV14}:  
\begin{align} \label{def:psi}
\left\{
\begin{array}{l} 
\partial_t (t\psi) = u_0^1 + \nu \partial_{yy}(t\psi) \\ 
\lim_{t \searrow 0} (t\psi) = 0.  
\end{array}
\right.
\end{align}
After applying this choice and using $U^1 = -(1+\partial_y \psi)(\partial_Y - t\partial_X) \Phi$ and $U^2 = \partial_X \Phi$, \eqref{eq:midstepOm} becomes 
\begin{align} \label{def:Omegnon}
\left\{
\begin{array}{l}
\partial_t \Omega + g\partial_Y \Omega + (1+\partial_y\psi) \grad^\perp \Phi_{\neq} \cdot \grad \Omega  = \nu \partial_{XX}\Omega + \nu(1+\partial_y \psi)^2(\partial_Y - t\partial_X)^2 \Omega \\ 
\Delta_t \Phi := \partial_{XX}\Phi + (1+\partial_y \psi)^2(\partial_Y - t\partial_X)^2 \Phi + \partial_{yy}\psi \partial_y \psi (\partial_Y - t\partial_X) \Phi = \Omega, \\ 
g = \partial_t \psi - \nu \partial_{yy}\psi = \frac{1}{t}(u_0^1 - \psi). 
\end{array}
\right.
\end{align}
Notice the crucial cancellation that eliminated powers of $t$ from the nonlinearity; such a nice structure is not present in 3D.   
The key difference between \eqref{def:Omegnon} and \eqref{pinguin}, is that the velocity field contains only $\Phi_{\neq}$, rather than $\Phi$: drift due to the slowly decaying shear flow has been removed from the equation for $\Omega$.
The presence of this perturbation is now being felt indirectly through $\psi$ via the forcing term \eqref{def:psi} (in particular, note that $\psi$ is an unknown that must be solved for along with $\Omega$).  
One can liken the coordinate transform \eqref{def:nonlin2Dcoord}, \eqref{def:psi} to a kind of gauge transformation.
Note that this coordinate change has made \eqref{def:Omegnon} significantly more nonlinear. 
In the proofs of \cite{BM13,BMV14}, further governing equations are derived for $\Psi'(t,Y) = \partial_y\psi(t,y)$ and $G(t,Y) = g(t,y)$ in order to obtain estimates on these nonlinear contributions. 

Notice that, like $\Delta_L$, $\Delta_t$ is not elliptic. Indeed, the symbol of $\Delta_t$ (as a pseudo-differential operator), 
\begin{align}
\sigma(\Delta_t)(Y,k,\eta) = -k^2 -(1+\partial_y\psi)^2(\eta-kt)^2 + i \partial_y\psi \partial_{yy}\psi (\eta-kt), \label{def:SymDeltat}
\end{align}
is degenerate at the frequency $\eta = kt$ -- the critical times.  
As alluded to above, the precise form of the coordinate transform \eqref{def:nonlin2Dcoord} was motivated by ensuring that the loss of ellipticity in \eqref{def:SymDeltat} still occurs at the same critical times as the Couette flow. 
This alone is not sufficient for us to consider $\Delta_t$ to be a small perturbation of $\Delta_L$; that requires $\partial_y \psi$ sufficiently small and a variety of arguments which carefully respect the precise way ellipticity is lost. 
Many variations of such elliptic estimates have appeared in varying levels of complexity in all of the works on Couette flow \cite{BM13,BMV14,BGM15I,BGM15II,BGM15III,BVW16}. 
Related arguments also arise when studying the enhanced dissipation due to the dissipative terms on the left-hand side of \eqref{def:Omegnon}. 

\bigskip

\noindent \textit{The 3D case}. In dimension three, the lift-up effect is so strong that the coordinate transform ends up being much larger and we need to account for it via a refined coordinate transform. 
In \cite{BGM15I,BGM15III}, the following ansatz is made in analogy with that made above:  
$$
\left\{
\begin{array}{l}
X = x - t y - t \psi(t,y,z) \\
Y = y + t \psi(t,y,z) \\
Z = z. 
\end{array}
\right.
$$
In \cite{BGM15II}, one also needs to make a more complicated transformation in $Z$ due to the very large size of the streak; see therein for details. 
Via a similar derivation as we applied above, one derives the following PDE for the evolution of $\psi$: 
$$
\frac{d}{dt}(t\psi) + U_0 \cdot \nabla (t\psi) = u_0^1 - t u_0^2 + \nu \Delta(t\psi),
$$
Notice the nonlinear transport and especially the lift-up effect seen through the presence of $-t u_0^2$. 
%These new coordinates allow to close energy estimates; however, the differential operators $\nabla$ and $\Delta$ become effectively quasilinear when expressed in $(X,Y,Z)$, leading to some further complications.

\subsection{Fourier multiplier norms} \label{sec:FourerMult}
One of the main technical tools employed in all of \cite{BM13,BMV14,BGM15I,BGM15II,BGM15III,BVW16}, as well as \cite{zillinger2016I,Zillinger2017,Bedrossian16} are various norms defined via carefully designed Fourier multipliers: for example, in 3D, defining $Q^i = \Delta_L U^i$
\begin{equation}
\label{normm}
\norm{A^i(t,\grad)Q^i}_{L^2} = \left( \sum_{k,\ell} \int \abs{A^i(t,k,\eta,\ell) \widehat{Q^i}(t,k,\eta,\ell)}^2 d\eta\right)^{1/2}.
\end{equation}
The multipliers are time-dependent, and typically encode two types of estimate: one is $L^\infty$ in time (global bound), the other $L^2$ in time (dissipation). This follows from the simple identity
\begin{align*}
\underbrace{\norm{A^i(t,\grad)Q^i(t) }_{L^2}^2}_{\displaystyle \mbox{uniform bound}} & =  \norm{A^i(0,\grad)Q^i(0)}_{L^2}^2 + 2 \underbrace{\int_0^t \langle \partial_t A^i(t,\grad)Q^i \, , \, A^i(t,\grad)Q^i \rangle\,ds}_{\displaystyle \mbox{dissipation term}} \\
 & + \underbrace{2 \int_0^t  \langle A^i(t,\grad)\partial_t Q^i\, , \, A^i(t,\grad) Q^i \rangle \,ds}_{\displaystyle \mbox{computed through the equation;  absorbed in various ways}}.
\end{align*}
Usually, the multipliers $A^i$ are designed as the product of several kinds of Fourier multipliers. 

%The kind of Fourier multiplier norm discussed in [REF] above is just one example of the kind of multipliers we have used in our works (and those used by Zillinger \cite{Zillinger2014} [REF].  
It turns out that a great diversity of Fourier multipliers are needed in the estimates of \cite{BM13,BMV14,BVW16,BGM15I,BGM15II,BGM15III,Zillinger2014} (and also the works in kinetic theory \cite{Bedrossian16,Bedrossian17}); each work uses a specific, different combination of them. 
 They generally can be classified into the following types: 

\begin{itemize} 
\item ``Ghost multipliers''. These multipliers are those which are uniformly bounded above and below in both $t$ and $\nu$; a Fourier-side analogue of Alinhac's ghost energy method \cite{Alinhac01}. These do not alter the norm topology and are hence relatively easy to apply. Variants have been applied throughout the works; they have been used in e.g. \cite{Zillinger2014,BGM15I,BGM15II,BGM15III,BVW16}.   

\item ``Nonlinear cascade multipliers''. These multipliers become weaker in time to allow the solution to lose large amounts of regularity due to a potential frequency cascade. Variants have been used in \cite{BM13,BMV14,BGM15I,BGM15II,Bedrossian16} and have each been designed based on weakly nonlinear toy models such as \eqref{ODETM}. They are much more difficult to use than ghost multipliers, especially those in \cite{BM13,BMV14,BGM15II}.

\item ``Steady loss/gain multipliers''. These multipliers steadily become weaker or stronger over long periods of time, for example, the standard Gevrey/analytic regularity multipliers such as $e^{\lambda(t)\abs{\grad}^s}$. However, the works \cite{BMV14,BGM15I,BGM15II,BGM15III} also use a variety of other multipliers which lose or gain in an anisotropic (in frequency and time) way.  %(for example, \cite{BGM15I,BGM15II} uses a multipilier which loses about $e^{c\abs{\partial_Y}^{1/2}}$ regularity steadily over the span of time $\abs{\partial_Y}^{1/2} \lesssim t \lesssim \abs{\partial_Y}$). 

\item ``Singular limit multipliers''. The most important multiplier in \cite{BGM15III} is different. This multiplier is uniformly bounded in $t$ and frequency, but \emph{not} uniformly in $\nu$. The multiplier was used to estimate the vortex stretching, which is unbounded as $\nu \rightarrow 0$ (but is bounded for all $\nu > 0$). 
We have elected to make this a separate classification of multiplier because we believe that such norms have a high probability of being useful in the future. 
%For each $\nu$, the viscosity eventually dominates and eliminates the enstrophy growth, however, the time-scale for this diverges as $\nu \rightarrow 0$ and the 3D Euler equations experience a linear-in-time, unbounded growth of enstrophy (this unbounded growth should be likened to the lift-up effect, but for enstrophy instead of energy). The multiplier employed in \cite{BGM15III} weakens the norm in a way which matches this transient growth. 
\end{itemize}

\subsubsection{An example of a ghost multiplier: enhanced dissipation}

Consider the simple linearized PDE: 
\begin{align}
\partial_t Q = \nu \Delta_L Q. \label{eq:toylin}
\end{align}
As we saw in Section~\ref{sec:ED2D}, this PDE is easily solved by taking the Fourier transform and integrating. 
However, once the nonlinearity is added back in, such an approach is not tenable. Indeed, semi-group methods will usually fail as $\nu \to 0$ since the problem becomes quasilinear, whereas energy estimates seem a much more promising tool.
Therefore, it is very natural to look for an energy method approach to getting enhanced dissipation estimates on \eqref{eq:toylin}. 
An approach using multipliers of the type we termed ``steady loss/gain'' was introduced in \cite{BMV14} and adapted in \cite{BGM15I,BGM15II}; however, while this approach easily yields very strong estimates, it is only suitable for proofs involving high regularity. 
An alternative approach, which seems to be much simpler and suitable for low regularity, was put forward in \cite{BGM15III} (and later applied in \cite{BVW16,Bedrossian17}). 
To use this approach, define the multiplier via the linear ODE for $k \neq 0$: 
%\begin{subequations} 
\begin{align*}
\frac{\dot{M}}{M} & = - \frac{\nu^{1/3}}{\left[ \nu^{1/3} |t-\frac{\eta}{k} | \right]^{2} + 1} \\ 
M(0,k,\eta) & = 1. 
\end{align*}
%\end{subequations} 
Notice that there is a constant $c$ (independent of $k$, $\eta$, $t$, and $\nu$) such that $c < M(t,k,\eta) \leq 1$, and hence this multiplier is of type we termed ``ghost multiplier''. 
In particular, its presence does not change a norm: 
\begin{align}
\norm{M(t,\grad) \brak{\grad}^\sigma f}_{L^2} \approx \norm{\brak{\grad}^\sigma f}_{L^2}. \label{ineq:Mequiv}
\end{align}
The crucial property that $M$ satisfies is: 
\begin{align}
1&\lesssim  \nu^{-1/6}\left(\sqrt{-\dot{M}M(t,k,\eta)}+\nu^{1/2}|k,\eta-kt|\right) \quad \mbox{for } k\neq0, 
\end{align}
which implies that 
\begin{align}
\norm{f_{\neq}}_{L^2}^2 \lesssim \nu^{-1/3}\left(\norm{ \sqrt{-\dot{M} M} f_{\neq}}_{L^2}^2 + \nu \norm{\grad_L f_{\neq}}_{L^2}^2\right). \label{ineq:Mprop}
\end{align}
Now consider the following simple energy estimate on \eqref{eq:toylin}: 
\begin{align}
\frac{1}{2}\frac{d}{dt}\norm{M(t,\grad) Q_{\neq}}_{L^2}^2 = -\norm{ \sqrt{-\dot{M} M} Q_{\neq}}_{L^2}^2 - \nu \norm{\grad_L M Q_{\neq}}_{L^2}^2, 
\end{align}
hence, 
\begin{align}
\frac{1}{2}\norm{M(t,\grad) Q_{\neq}(T)}_{L^2}^2 + \int_0^T \norm{\sqrt{-\dot{M}M} Q_{\neq}}_{L^2}^2 dt + \nu \int_0^T \norm{\grad_L M Q_{\neq}}_{L^2}^2 dt = \frac{1}{2}\norm{M(0,\grad) Q_{\neq}(0)}_{L^2}^2. 
\end{align}
Applying \eqref{ineq:Mequiv} and \eqref{ineq:Mprop} gives us the decay estimate 
\begin{align}
\int_0^\infty \norm{M Q_{\neq}(t)}_{L^2}^2 dt \lesssim \nu^{-1/3} \norm{Q(0)}_{L^2}^2.  
\end{align}
Hence, this estimate scales with the correct $\nu^{-1/3}$ characteristic time-scale observed using the Fourier approach. 
Further, one can adapt this method to obtain exponential decay rates such as $e^{-\delta\nu^{1/3}t}$ for sufficiently small $\delta$ \cite{Bedrossian17}.

\subsubsection{An example of a nonlinear cascade multiplier: decay vs regularity in the simple 2D case} 
We build up here on the toy model derived in Section~\ref{TM} to derive a multiplier which will control the frequency cascade scenario.
We focus on the 2D case since the 3D case is significantly more technical. 
Recall that in \S\ref{TM}, we came up with the following simplified model
\begin{align} \label{def:OmegnonToy}
\left\{
\begin{array}{l}
\partial_t \Omega + \grad^\perp \Phi_{H} \cdot \grad \Omega_L  = \nu \Delta_L \Omega \\ 
\Delta_L \Phi = \Omega. 
\end{array}
\right.
\end{align}
By repeating calculations from \S\ref{IDOM2D}, we have the following decay estimate for the effective velocity field due to inviscid damping: 
\begin{align}
\norm{\grad^\perp \phi_{\neq}}_{H^s} \lesssim \frac{\norm{\Omega(t)}_{H^{s+3}}}{\brak{t}^2}. \label{ineq:regloss2d}
\end{align}
It follows that if we have a uniform-in-time bound on $\Omega$ in a suitable $H^s$ space, then the nonlinearity decays rapidly and the nonlinear dynamics match the linear dynamics predicted by \eqref{linOmega2D}. The metastability discussed in \S\ref{sec:shear2Dstab} will also hold accordingly (at least after we quantify the enhanced dissipation in a suitable sense). 
However, due to the regularity loss in \eqref{ineq:regloss2d}, it is not clear how to obtain this uniform bound -- the loss is far too large even to get a uniform bound in analytic regularity via a Cauchy-Kovalevskaya style argument (even to get $t^{-1}$ decay requires the loss of two derivatives whereas a Cauchy-Kovalevskaya argument could handle at most the loss of one).

We see that the estimate \eqref{ineq:regloss2d} is too simplistic to accomplish anything (though it is optimal if one restricts oneself to standard $H^s$ norms and pointwise-in-time estimates).
Let us return to \eqref{ODETM} and try to more precisely estimate the kind of regularity loss that could occur. 
 As we are only interested in an upper bound, we will drop the signs in \eqref{ODETM}. If we assume that $t \sim \eta/k$, then \eqref{ODETM} simplifies further to the following, assuming that $\Omega_L = O(\kappa)$: 
\begin{align*}
%\partial_t \widehat{\Omega}_H(t,k,\eta) & = -\frac{\eta}{(k-1)^2 + (\eta-(k-1)t)^2} \widehat{\Omega}_H(t,k-1,\eta) \Omega_L(1,0) \\ 
% \partial_t \widehat{\Omega}_H(t,k-1,\eta) & = \frac{\eta}{k^2 + (\eta-kt)^2} \widehat{\Omega}_H(t,k,\eta) \Omega_L(-1,0).
\partial_t \widehat{\Omega}_H(t,k,\eta) & = \kappa\frac{\eta}{(k-1)^2 + (\eta-(k-1)t)^2} \widehat{\Omega}_H(t,k-1,\eta) \\ 
\partial_t \widehat{\Omega}_H(t,k-1,\eta) & = \kappa\frac{\eta}{k^2 + (\eta-kt)^2} \widehat{\Omega}_H(t,k,\eta). 
\end{align*} 
We will further assume that $\eta \gtrsim k^2$ (so that the coefficient in the second equation is large) and that $\abs{t - \eta/k} \lesssim \frac{\eta}{k^2}$, which implies 
that $t$ is near the critical time $\eta/k$ and away from the critical time $\eta/(k-1)$. 
Applying these simplifications derives the \emph{toy model} of \cite{BM13}: 
\begin{subequations} \label{2Dtoy}
\begin{align}
%\partial_t \widehat{\Omega}_H(t,k,\eta) & = -\frac{\eta}{(k-1)^2 + (\eta-(k-1)t)^2} \widehat{\Omega}_H(t,k-1,\eta) \Omega_L(1,0) \\ 
% \partial_t \widehat{\Omega}_H(t,k-1,\eta) & = \frac{\eta}{k^2 + (\eta-kt)^2} \widehat{\Omega}_H(t,k,\eta) \Omega_L(-1,0).
\partial_t \widehat{\Omega}_H(t,k,\eta) & = \kappa\frac{k^2}{\eta} \widehat{\Omega}_H(t,k-1,\eta) \\ 
\partial_t \widehat{\Omega}_H(t,k-1,\eta) & = \kappa\frac{\eta}{k^2 + (\eta-kt)^2} \widehat{\Omega}_H(t,k,\eta). 
\end{align}
\end{subequations}
It is possible to find an approximate supersolution of the coupled ODE system~\eqref{2Dtoy}, which, after iterating over all critical times, quantifies precisely how bad the echo cascade can be \cite{BM13}. Next to a critical time $t \sim \frac{\eta}{k}$, the supersolution $w_C$ controls the critical frequency $k$, and the supersolution $w_{NC}$ the non-critical frequency $k-1$:
\begin{align*}
\abs{\Omega_H(t,k,\eta)} & \lesssim w_C(t,k,\eta) \\ 
\abs{\Omega_H(t,k-1,\eta)} & \lesssim w_{NC}(t,k,\eta), 
\end{align*}
where $w_C$ and $w_{NC}$ are functions which satisfy the following (near $t \approx \eta/k$),  
\begin{align*}
w_C\left(\frac{\eta}{k}+\frac{\eta}{k^2},k,\eta\right) & \approx w_{NC}\left(\frac{\eta}{k}+\frac{\eta}{k^2},k,\eta\right) \approx \left(\frac{\eta}{k^2}\right)^{1 + O(\kappa)}w_{NC}\left(\frac{\eta}{k} - \frac{\eta}{k^2}k,\eta\right) 
\end{align*}
and 
\begin{align}
\frac{w_{NC}(t,k,\eta)}{w_{C}(t,k,\eta)} \approx \frac{\eta}{k^2(1 + \abs{t - \frac{\eta}{k}})}.  \label{ineq:toyimbal}
\end{align}
The first approximate identity shows that the supersolution estimates that the total growth of the two modes is roughly comparable. 
That is, both the critical and non-critical frequencies grow by a factor of about $O(\frac{\eta}{k^2})^c$ for some fixed $c$. 
Iterating over all critical times satisfying $\eta \gtrsim k^2$ predicts a growth in the mode $(1,\eta)$ by time $t\geq \eta$ like the following (applying Stirling's formula; see \cite{BM13} for more details):
\begin{align*}
\left(\frac{\eta^{\sqrt{\eta}}}{(\sqrt{\eta}!)^2}\right)^c \sim \eta^{-c/2} e^{2c\sqrt{\eta}}.
\end{align*}
Hence, \eqref{2Dtoy} precisely predicts that the amount of enstrophy at frequencies comparable to $\eta$ can be amplified by roughly $e^{2c\sqrt{\eta}}$ -- this large amplification suggests that Gevrey-2 regularity might be necessary in order to control the nonlinear effects for long times (if one does not have viscosity). 
This is the origin of the high regularity requirement in \cite{BM13,BMV14}. 

The approximate identity \eqref{ineq:toyimbal} on the other hand shows that, very close to $\eta \approx tk$, the non-critical frequencies could potentially be much larger than the critical frequencies. 
This regularity imbalance is only possible due to the fact that the critical frequency does not interact directly with itself in the leading order dynamics predicted by \eqref{2Dtoy}. 
This special non-interaction can be considered a \emph{null form} which is absolutely crucial to the proof of \cite{BM13,BMV14} (it arises from the $\grad^\perp \Phi \cdot \grad$ structure in the original equations). 
Such null forms that dictate the precise way information moves between different frequencies also prove to be essential to the 3D works, though there the structures are more complicated \cite{BGM15I,BGM15II,BGM15III}. 

The proof of \cite{BM13} is based on a weighted norm roughly of the form $A = e^{\lambda(t)\abs{\grad}^{s}}\frac{1}{w(t,k,\eta)}$ with $s \in (0,1)$ and $w$ designed to resemble the above supersolution near the critical times $t \sim \eta/j$ for any $j$ with $\abs{j}^2 \lesssim \abs{\eta}$ (after accounting for the subsequent multiplicative amplification through each critical time). 
That is, for $k = j$, $w \sim w_{C}$ and $k \neq j$, we take $w \sim w_{NC}$. 
The motivation for this choice can be interpreted as follows: $w(t,k,\eta)$ roughly predicts the ``worst-case'' transfer of information from low-to-high frequencies, and, if the enstrophy transfer is of roughly the same type as that predicted, we can expect that 
\begin{align*}
\abs{\Omega(t,k,\eta)} \lesssim w(t,k,\eta), 
\end{align*}
which ensures that the norm gets weaker at precisely the rate necessary to ensure $A\Omega$ remains uniformly bounded.

In the 3D works, the analogues of the toy model are far more complicated \cite{BGM15I,BGM15II}, however, these can be upper bounded by multipliers which are similar to the above $w$.
In the work on Vlasov-Poisson, a multiplier analogous to the above $w$ is used, however, it is also a singular limit (in the size of the data $\eps$), and so it is in some ways more complicated but in other ways simpler, as the resonances in Vlasov-Poisson have a simpler structure \cite{Bedrossian16}.

\subsection{Paraproducts and Gevrey regularity} 
\label{paraprod}
The paraproduct, introduced by Bony in \cite{Bony81},  
%\comment{(Pierre) We should add the reference J.-M. Bony, Calcul symbolique et propagation des singularit\'es pour les \'equations aux d\'eriv\'ees partielles non lin\'eaires, Ann. Sci. Ecole Norm. Sup. 14 (1981), 209--246.}
has been a ubiquitous tool in the analysis of the stability of the Couette flow.
In order to explain its principle, consider functions $f$ and $g$ of a real variable, and recall that the Fourier transform of their product is given by the convolution of their Fourier transforms:
$$
\widehat{fg}(\xi) = \int \widehat{f}(\xi-\eta) \widehat{g}(\eta) \, d\eta.
$$
The idea is now to use a smooth cut off function $\chi$ such that
$$
\chi(\xi,\eta) = 
\left\{ \begin{array}{ll}
1 & \mbox{if $|\eta| >> |\xi-\eta|$} \\
0 & \mbox{if $|\xi-\eta| >> |\eta|$}
\end{array} \right.
$$
in order to split the product
$$
\widehat{fg}(\xi) = \int \chi(\xi,\eta) \widehat{f}(\xi-\eta) \widehat{g}(\eta) \, d\eta + \int [1-\chi(\xi,\eta)] \widehat{f}(\xi-\eta) \widehat{g}(\eta) \, d\eta.
$$
In the right-hand side, the first summand corresponds to interactions between high frequencies of $g$ and low frequencies of $h$. This explains the following shorthand for the above decomposition:
$$
fg = f_{\operatorname{Lo}} \, g_{\operatorname{Hi}} + f_{\operatorname{Hi}} \, g_{\operatorname{Lo}}.
$$
The decomposition of all products in this way turns out to be very convenient, in that it allows to use different estimates for different interactions. Consider for instance the paraproduct decomposition of the convection term
$$
u \cdot \nabla u = u_{\operatorname{Lo}} \cdot (\nabla u)_{\operatorname{Hi}} + u_{\operatorname{Hi}} \cdot (\nabla u)_{\operatorname{Lo}}
$$
The second term should be better behaved, since the derivative hits low frequencies. This can be quantified through the following estimate
$$
\left\| u_{\operatorname{Hi}} \cdot (\nabla u)_{\operatorname{Lo}} \right\|_{H^s} \lesssim \| u \|_{H^s}^2.
$$
As a consequence, this particular term does not ``lose derivatives'', as opposed to the low-high interaction.

Basic paraproduct methods provide a straightforward proof of the product rule for fractional Sobolev spaces (see e.g. the appendix of \cite{TaoTextbook}),
and hence it is reasonable to expect that they also provide a useful tool for deducing product rules in Gevrey class.
Indeed, one has the following very useful product rule which emphasizes one of the very convenient properties of Gevrey regularity.
Specifically, only one of the two factors is measured in the highest norm in each of the two terms (analogous to the situation in Sobolev spaces).
We remark that this is false in analytic regularity $s=1$. 
\begin{thm}
Let $s \in (0,1)$, $p > d/2$, and $\lambda > 0$. There exists a $c= c(s) \in (0,1)$ depending only on $s$, such that
\begin{align}
\norm{fg}_{\cG^{\lambda,s}} \lesssim_{s,p} \norm{f}_{\cG^{\lambda,s}} \norm{\brak{\grad}^{p}g}_{\cG^{c\lambda,s}} + \norm{g}_{\cG^{\lambda,s}} \norm{\brak{\grad}^{p}f}_{\cG^{c\lambda,s}}.
\end{align}
\end{thm}
This theorem can be proved using paraproducts together with elementary inequalities of the general type $\abs{x-y}^s < c(s)\left(\abs{x}^s + \abs{y}^s\right)$ (see e.g.~\cite{BM13} for a proof).
We remark that when working in Gevrey class, it is necessary to split into three regions of frequency rather than two: $\abs{\xi} \ll \abs{\eta-\xi}$, $\abs{\eta-\xi} \ll \abs{\xi}$, and $\abs{\eta-\xi} \approx \abs{\xi}$.

\section{Related questions}

\subsection{Scattering and weak turbulence} \textit{Scattering} is a typical behavior for (defocusing) nonlinear dispersive equations set in the whole space: global solutions look asymptotically linear. To be more specific, consider for instance the nonlinear Schr\"odinger equation
\begin{equation}
\label{NLS} \tag{NLS}
\left\{
\begin{array}{l}
i \partial_t \psi - \Delta \psi = - |\psi|^2 \psi \\
\psi(t=0) = \psi_0 
\end{array}
\right.
\end{equation}
set on the whole space $\mathbb{R}^3$: $\psi = \psi(t,x) \in \mathbb{C}$, with $(t,x) \in \mathbb{R} \times \mathbb{R}^3$. As $t \to \pm \infty$, $\psi(t)$ converges locally to zero, and looks increasingly like a linear solution:
$$
\left\| \psi(t) - e^{it\Delta} f_{\pm} \right\|_{H^1} \rightarrow 0 \quad \mbox{as $t \to \infty$}
$$
(in general, $f_{+}$ and $f_{-}$ are different). 

\medskip

Recall that, in dimension $2$, for $\nu = 0$, and for small Gevrey data,
\begin{align}
\omega(t,x,y) \approx \omega_{\pm\infty}(t,x-ty-tu_\infty(y),y) \quad\quad \mbox{as $t \to \pm \infty$}.
\end{align}
The similarities with scattering for~\eqref{NLS} are striking: in both cases, pointwise decay ($\psi \to 0$ and $u^2 \to 0$ respectively); the system relaxes as $t \to \pm \infty$ (to a linear solution), though no dissipation appears in the system, which is fully time-reversible; and the scattering states $\omega_{\pm \infty}$ differ in general.

The mechanism of relaxation for~\eqref{NLS} is the spreading of waves, which causes them to decay, and ultimately suppresses nonlinear interactions. No such mechanism is available for~\eqref{E}, since the domain is compact in $X$ (and very little spreading in $Y$ occurs). However, the energy is spreading~\textit{in frequency}: drift towards high frequencies (and the Biot-Savart law) is the cause of inviscid damping.

This is very reminiscent of a phenomenon called~\textit{weak turbulence}. It is expected for~\eqref{NLS} set on a compact domain, and consists of a drift of energy towards high frequencies as $t\to \infty$, causing the system to relax. However, no mathematical evidence is known for this phenomenon, which, furthermore, should be valid in a statistical sense only.

\subsection{Landau damping} \label{sec:LandauMixing} 
The Vlasov-Poisson equation is a kinetic model describing the evolution of a density of particles $f(t,x,v)$, which, at time $t$, have position $x$ (usually taken $x \in \Torus^d_x$ or $x \in \Real^d_x$) and velocity $v \in \Real^d_v$:
\begin{equation}
\label{VP}
\tag{VP}
\left\{
\begin{array}{l}
\partial_t f + v \cdot \nabla_x f + E \cdot \nabla_v f = 0 \\
E(t,x) = \nabla_x \Delta^{-1}_x ( \int_{\Real^d} f dv - \overline{\rho}) \\
f(t=0) = f_{in}.
\end{array}
\right.
\end{equation}
where $\overline{\rho}$ is a neutralizing background density (usually fixed to be a constant) which models the presence of the ions in the plasma \cite{BoydSanderson}. 
At least to start, in physical applications, one studies disturbances of spatially homogeneous background $f(t,x,v) = f^0(v) + h(t,x,v)$
with the mean-zero condition $\int_{\Torus^d} \int_{\Real^d} h(t,x,v) dv dx = 0$ and $\overline{\rho} = \int f^0(v) dv$. 
The equations hence become
\begin{equation}
\label{VP2}
\left\{
\begin{array}{l}
\partial_t h + v \cdot \nabla_x h + E \cdot \left(\nabla_v f^0 + \grad_v h\right) = 0 \\
E = \nabla_x \Delta^{-1}_x \int_{\Real^d} h dv \\  
h(t=0) = h_{in}.
\end{array}
\right.
\end{equation}
If $h$ solved the free transport equation with smooth data: 
$$
\partial_t h + v \cdot \nabla_x h = 0,
$$
then it is easy to see that $E \to 0$ as $t \to \infty$ (it suffices to take the Fourier transform and proceed exactly as in the Couette flow). 
This phenomenon is called \textit{Landau damping}; just like inviscid damping, it gives decay in the absence of a dissipative mechanism. 
We also remark that it is a variant of velocity averaging in a slightly different guise; see e.g. \cite{GolseEtAl1985,GolseEtAl1988}. 
In 1946, Landau proved that the Landau damping holds for the linearized Vlasov equations on $x \in \Torus^d$ for $f^0$ a Maxwellian \cite{Landau46} and Penrose \cite{Penrose} later extended this to a wide class of other backgrounds. See e.g. \cite{VKampen55,Degond86,Glassey94,Glassey95} for other studies on the linearized Vlasov equations. We remark that it was apparently Van Kampen \cite{VKampen55} who first pointed out that Landau damping is related to mixing in phase space. It is one of the most fundamental effects in the kinetic theory of plasmas \cite{Ryutov99,BoydSanderson} and is also potentially relevant to stellar mechanics \cite{LyndenBell67}.

It was first proved that there exists solutions to the nonlinear Vlasov equations \eqref{VP2} on $x \in \Torus$ which display Landau damping by Caglioti and Maffei \cite{CagliotiMaffei98}, extended further later in \cite{HwangVelazquez09}. 
A major breakthrough came in Mouhot and Villani's work \cite{MouhotVillani11}, which showed that  on $x \in \Torus^d$,  the dynamics of all sufficiently small and (sufficiently smooth) Gevrey class solutions matched that of the linearized Vlasov equations for all $t$ (and hence, all such solutions display Landau damping at essentially the same rate as that predicted by Landau and Penrose's work). 
The proof was later simplified and extended to a wider range of Gevrey regularity by Mouhot and two of the authors in \cite{BMM13}, which was later adapted to the relativistic plasma case \cite{Young14}.
In \cite{MR3437866,FGVG} Landau damping in Sobolev regularity was proved for variants of the Vlasov-Poisson model which have much weaker resonances (see also \cite{MR3471147}). 
Later, dispersion and phase mixing were combined to prove Landau damping on $\Real^3_x \times \Real^3_v$ in \emph{Sobolev} regularity for Vlasov-Poisson in \cite{BMM16}. In particular, it was proved that dispersive effects are able to limit the effect of the plasma echo resonances.
However, in \cite{Bedrossian16}, it was proved that on $\Torus_x \times \Real_v$ there exists solutions arbitrarily close to homogeneous equilibrium in Sobolev spaces which deviate arbitrarily far from the linearized Vlasov predictions due to long chains of plasma echo resonances.
Recently in \cite{Bedrossian17}, it was shown that nonlinear collisions can suppress these resonances in a manner analogous to the results of \cite{BGM15III,BVW16}. 

\section{Open problems}

The results of \cite{BM13,BMV14,BGM15I,BGM15II,BGM15III,BVW16} study what is essentially the easiest problem in hydrodynamic stability at high Reynolds in the absence of unstable spectra in the linearization.
However, now that this case has been studied, it makes sense to ask what harder problems can be considered. 
Moreover, despite the rather in-depth studies, there are still some unanswered questions even on Couette flow.

\subsection{Nonlinear instabilities and transition thresholds: open questions on Couette}

One of the most important directions to explore further is the role of nonlinear instabilities in subcritical transition, and moreover, how the regularity of initial data may or may not affect the dynamics. In particular, an important problem regarding Couette flow would be to prove the sharpness of the transition thresholds estimated in 2D \cite{BVW16} and 3D \cite{BGM15III} for Sobolev regularity data.
A related problem would be studying nonlinear resonances in the inviscid 2D Euler equations near Couette to provide a hydrodynamic analogue of the study of plasma echoes in \cite{Bedrossian16}.
Similarly, one could evaluate carefully how things change in still lower regularity, for example, see \cite{LinZeng11} for 2D Euler.  
At least for the case of the 3D Couette flow, there is a variety of quantitative and anecdotal evidence to suggest that, indeed, regularity does play an important role in the subcritical transition. The numerical studies of Reddy et. al. \cite{ReddySchmidEtAl98} suggest specifically that the threshold in \cite{BGM15III} is sharp, however, numerical studies of these problems are very difficult and prone to overestimate the transition thresholds so it is difficult to say confidently either way. 
It makes sense to try to connect any such analysis to existing formal asymptotic works such as \cite{Vanneste02,VMW98} and physical experiments \cite{SchmidHenningson2001,Yaglom12}. 
Finally, one would naturally wish to understand better the `secondary instabilities' of streaks and the role played in the next step in the transition from nearly linear to nonlinear dynamics. 

\subsection{Influence of boundaries and low wavenumbers} \label{sec:openboundariesx}

The physical experiments most closely approximating the Couette flow setting are in bounded channels $y \in [-1,1]$ and long channels better modeled by $x \in \Real$. 
The boundary conditions at the top and bottom are naturally $u = (1,0,0)$ at $y=1$ and $u = (-1,0,0)$ at $y=-1$ after proper normalization (we will refer to this case as the infinite channel problem).  
It has been observed that the presence of physical boundaries could potentially introduce a variety of new kinds of nonlinear instabilities (see \cite{SchmidHenningson2001,Yaglom12} and the references therein). 
Moreover, the presence of long wave-numbers in $x$ presents another challenge. 
It has been proved that many \emph{planar} shear flows in the channel that are not Couette flow suffer from long-wave instabilities at high Reynolds number
(\cite{Grenier2016,FriedlanderEtAl2006} and the references therein)
Understanding both (or either) the long wave effects and the boundary effects in 3D is an important goal moving forward. 
Even at the level of 2D linearized Euler in a channel, it has been observed that boundaries can add new complications \cite{zillinger2016I,WeiZhangZhao15}. 

\begin{remark} 
Studying nonlinear stability for 2D Navier-Stokes at high Reynolds number or 2D Euler in a bounded channel (e.g. $\Torus \times [-L,L]$ or $\Real \times [-L,L]$) may not be of direct physical relevance. 
As discussed above, real 3D shear flows undergo strong 3D instabilities at high Reynolds number -- the 2D dynamics are rarely observed. 
Applications of 2D fluid mechanics to plasma physics (via `drift-kinetic' or `gyro-kinetic' scenarios) or atmospheric dynamics would generally not involve such boundaries, or, would  have  boundary conditions (and boundary layers) which are different from the usual no-slip and no-penetration conditions. 
Nevertheless, this case is very interesting and important to study from a mathematical viewpoint, as one can use this case as an intermediate result on the path to understanding the case of 3D boundaries (which is of high physical interest).   
\end{remark} 

\subsection{Stability and subcritical instability for more general problems} 

A variety of other settings in fluid mechanics are of more direct physical relevance than the Couette flow.
In two dimensions, the most important configuration to study are vortices in 2D Navier-Stokes and Euler (that is, radially symmetric configurations of vorticity). 
In analogy with shear flows, for the Euler equations all such configurations are equilibria whereas for Navier-Stokes all such solutions reduce to the heat equation. 
See \cite{Gallay2017,LiWeiZhang2017,BCZV17} and references therein for progress. 
In three dimensions there are several configurations of major relevance.
Perhaps the most famous of the problems is the pipe flow studied by Reynolds in his original experiments \cite{Reynolds83}: that is the equilibrium corresponding to pressure-driven flow in a cylindrical pipe (preferably an infinite pipe, but the periodic model would likely be studied first).
Other good examples are wall bounded shear flows (e.g. semi-infinite domains with a parabolic shear flow above a flat plate), Couette in an infinite channel (bounded in $y$), or other planar shear flows in an infinite channel \cite{DrazinReid81}. 
Naturally one could also study vortex columns. 
For each of these problems, one can formulate analogues of Questions \ref{question} and \ref{question2} in various regularities and pursue questions similar to those we studied in the Couette flow.  
Similarly, various kinds of subcritical instabilities will arise in these problems and it would be important to isolate them and study them in more detail. 
Moreover, many similar problems exist in kinetic theory and magneto-hydrodynamics problems, for example, understanding Landau damping in galactic dynamics or in plasmas in the presence of magnetic fields. 

As we saw above, the stability of the Couette flow has been analyzed at depth (though some questions remain).
One of the important reasons why this was possible is that the linearized problem can be solved explicitly.
For more general problems this is no longer possible, and indeed, even the resulting \emph{linear} problems are extremely difficult to analyze in detail.
Linear works studying enhanced dissipation and inviscid damping in 2D Navier-Stokes and Euler have been quite technical and the theory is still in its infancy (see e.g. \cite{zillinger2016I,Zillinger2017,BCZV17,WeiZhangZhao15,WeiZhangZhao2017,Gallay2017,LiWeiZhang2017,WZZK2017,IYM17}). 
We are not aware of any analogous works for 3D Navier-Stokes or 3D Euler, indeed, this kind of detailed information about the linearized operators in 3D seems to be completely open (except for Couette flow) --  even for the original problem of Reynolds \cite{Reynolds83,DrazinReid81} after 120 years of research. 
Moreover, our works on Couette flow show that nonlinear results do not easily follow from linear results. 
However, it is our hope that the ideas and methods which were put forward in the analysis of the Couette flow can lead to a much more general understanding and eventually lay the groundwork for a wider theory. 

\subsection{Large-time behavior of Euler for $d=2$} The existence of global solutions of the Euler equations~\eqref{E} in dimension $d=2$ is known under weak assumptions on the data. However, very little is known qualitatively about these solutions.
Until recently, the only known solutions whose asymptotic behavior could be established were stationary or explicit solutions! 
Gevrey solutions close to Couette are the first nontrivial solutions whose behavior can be described precisely \cite{BM13}. 
The works of e.g. \cite{KiselevSverak2014,Nadirashvili1991} and related norm growth results also provide glimpses into interesting long-term dynamics, though significantly less information is obtained on solutions in these regimes. 
The work of e.g. \cite{CastroEtAl2016} moreover construct a variety of interesting smooth, rigidly rotating solutions. Any solutions of these types clearly do not involve vorticity mixing as $t\rightarrow \infty$, and so we can expect a large set of solutions which do not display any mixing. 
On the other hand, results of \cite{BM13} show that \emph{all} sufficiently smooth solutions in the vicinity of Couette flow which are not shear flows experience some vorticity mixing as $t \rightarrow \infty$ and in particular, are \emph{not} precompact in $L^2(\Torus \times \Real)$. 
In light of \cite{BM13}, statistical mechanics considerations, and a variety of numerical simulations, it is potentially reasonable to make the following, somewhat vague, conjecture which suggests that this situation is generic in a suitable sense.
This conjecture was first communicated to us by Vladimir Sverak. 

\begin{conjecture} 
The `generic' solution to 2D Euler in vorticity form on $\Torus^2$ is such that the orbit $\set{\omega(t)}_{t \in \Real}$ is not precompact in $L^2(\Torus^2)$.
Here `generic' could be interpreted in the sense of Baire category or in an appropriate probabilistic sense, such as randomizing initial data in a suitable manner. 
\end{conjecture}

\bibliographystyle{abbrv} \bibliography{eulereqns,eulereqns_vlad}

\end{document}